\documentclass[journal]{IEEEtran}
\ifCLASSINFOpdf
\else
\fi
\usepackage[utf8]{inputenc}
\usepackage{graphicx}
\usepackage{amsmath} 
\usepackage{mathtools} 
\usepackage{amssymb}  
\usepackage{psfrag}
\usepackage{mathrsfs}
\usepackage{color}
\usepackage{cases}
\usepackage{empheq}
\usepackage{pgfplots}
\usepackage{tikz}
\usepackage{pgf-extrect}
\usepackage[authormarkup=none]{changes}
\usetikzlibrary{calc}
\usetikzlibrary{intersections, positioning,fit,arrows.meta,backgrounds}
\usepackage[final,activate={true,nocompatibility},shrink=10,stretch=10]{microtype}
\usepackage{bm} 

\newtheorem{thm}{Theorem}

\newtheorem{rem}{Remark}

\renewcommand{\d}{\mathrm{d}}
\newcommand{\e}{\mathrm{e}}
\newcommand{\prho}{\bm{\rho}}
\newcommand{\psigma}{\bm{\sigma}}
\newcommand{\pG}{\bm{G}}
\newcommand{\pH}{\bm{H}}
\newcommand{\pK}{\pmb{K}} 
\newcommand{\pa}{\bm{a}}
\newcommand{\pA}{\bm{A}}
\newcommand{\pxi}{\bm{\xi}}
\newcommand{\peta}{\bm{\eta}}
\newcommand{\ps}{\bm{s}}
\newcommand{\pz}{\bm{z}}
\newcommand{\pzeta}{\bm{\zeta}}
\newcommand{\dz}{\partial_{z}}
\newcommand{\dzeta}{\partial_{\zeta}}
\newcommand{\ie}{i.\,e., }

\newcommand{\wrt}{w.\,r.\,t. }
\newcommand{\ozeta}{\bar{\zeta}}
\newcommand{\oeta}{\bar{\eta}}

\newcommand{\oxi}{\bar{\xi}}
\newcommand{\orho}{\bar{\rho}}

\newcommand{\rightcond}[1]{\scriptsize\setlength{\arraycolsep}{0.5pt}\begin{array}[t]{rl}
#1
\end{array}}

\definechangesauthor[name=Simon,color=blue]{SK} 

\newcommand{\mathpushright}[2][1]{
	\mathmakebox[#1\displaywidth][r]{#2}
}
\newcommand{\mathpushleft}[2][1]{
	\mathmakebox[#1\displaywidth][l]{#2}
}

\newcommand{\tikzmark}[3][3pt]{\tikz[remember picture,baseline=(#2.base)] 
	\node[inner sep=#1,outer sep=0pt] (#2) {#3};}

\makeatletter
	\def\parsenode[#1]#2\pgf@nil{%
	    \tikzset{label node/.style={#1}}
	    \def\nodetext{#2}
	}
	
	\tikzset{
	    add node at x/.style 2 args={
	        name path global=plot line,
	        /pgfplots/execute at end plot visualization/.append={
	                \begingroup
	                \@ifnextchar[{\parsenode}{\parsenode[]}#2\pgf@nil
	            \path [name path global = position line #1-1]
	                ({axis cs:#1,0}|-{rel axis cs:0,0}) --
	                ({axis cs:#1,0}|-{rel axis cs:0,1});
	            \path [xshift=1pt, name path global = position line #1-2]
	                ({axis cs:#1,0}|-{rel axis cs:0,0}) --
	                ({axis cs:#1,0}|-{rel axis cs:0,1});
	            \path [
	                name intersections={
	                    of={plot line and position line #1-1},
	                    name=left intersection
	                },
	                name intersections={
	                    of={plot line and position line #1-2},
	                    name=right intersection
	                },
	                label node/.append style={pos=1}
	            ] (left intersection-1) -- (right intersection-1)
	            node [label node]{\nodetext};
	            \endgroup
	        }
	    },
	    add node at y/.style 2 args={
	        name path global=plot line,
	        /pgfplots/execute at end plot visualization/.append={
	                \begingroup
	                \@ifnextchar[{\parsenode}{\parsenode[]}#2\pgf@nil
	            \path [name path global = position line #1-1]
	                ({axis cs:0,#1}-|{rel axis cs:0,0}) --
	                ({axis cs:0,#1}-|{rel axis cs:1,1});
	            \path [yshift=1pt, name path global = position line #1-2]
	                ({axis cs:0,#1}-|{rel axis cs:0,0}) --
	                ({axis cs:0,#1}-|{rel axis cs:1,1});
	            \path [
	                name intersections={
	                    of={plot line and position line #1-1},
	                    name=left intersection
	                },
	                name intersections={
	                    of={plot line and position line #1-2},
	                    name=right intersection
	                },
	                label node/.append style={pos=1}
	            ] (left intersection-1) -- (right intersection-1)
	            node [label node] {\nodetext};
	            \endgroup
	        }
	    }
	}
\makeatother

\begin{document}
%
\title{Backstepping Control of Coupled Linear Parabolic PIDEs with Spatially-Varying Coefficients}
%
%
%

\author{Joachim~Deutscher and Simon~Kerschbaum
\thanks{J. Deutscher and S. Kerschbaum are with the Lehrstuhl f\"ur Regelungstechnik, Universit\"at
    Erlangen-N\"urnberg, Cauerstra\ss e 7, D-91058~Erlangen, Germany.
    (e-mail: \{joachim.deutscher,simon.kerschbaum\}@fau.de)}}

%
%

\markboth{IEEE TRANSACTIONS ON AUTOMATIC CONTROL,~Vol.~XX, No.~Y, JULY~2017}%
{Shell \MakeLowercase{\textit{et al.}}: Bare Demo of IEEEtran.cls for Journals}
%



\maketitle

\begin{abstract}
This paper considers the backstepping design of state feedback controllers for coupled linear parabolic partial integro-differential equations (PIDEs) of Volterra-type with distinct diffusion coefficients, spatially-varying parameters and mixed boundary conditions. 
The corresponding target system is  a cascade of parabolic PDEs with local couplings allowing a direct specification of the closed-loop stability margin. 
The determination of  the state feedback controller leads to kernel equations, which are a system of coupled linear second-order hyperbolic PIDEs with spatially-varying coefficients and rather unusual boundary conditions.
By extending the method of successive approximations for the scalar case to the considered system class, the well-posedness of these kernel equations is verified by providing a constructive solution procedure. This results in a systematic method for the backstepping control of coupled parabolic PIDEs as well as PDEs. The applicability of the new backstepping design method is confirmed by the stabilization of two coupled parabolic PIDEs with Dirichlet/Robin unactuated boundaries and a coupled Neumann/Dirichlet actuation.
\end{abstract}

\begin{IEEEkeywords}
Distributed-parameter systems, parabolic systems, coupled PIDEs, boundary control, backstepping.
\end{IEEEkeywords}

\section{Introduction}
\IEEEPARstart{I}{n} the last decade the \emph{backstepping approach} has emerged as a very powerful tool to solve stabilization problems for boundary controlled distributed-parameter systems (DPS) (see \cite{Kr08} and \cite{Me13} for an overview). The basic idea of this method is to utilize an invertible Volterra-type integral transformation for facilitating the controller design. 
As far as parabolic systems are concerned, the backstepping method for the state feedback design was first introduced in \cite{Liu03,Sm04} for scalar parabolic PDEs and PIDEs. Subsequently, these results were extended to parabolic systems with space- and time-dependent coefficients in \cite{Sm05,Me09}, to parabolic PDEs with Volterra nonlinearities in \cite{Vaz08,Vaz08a} and to higher-dimensional spatial domains in \cite{Me13,Jad15}. 

Recently, the backstepping control of \emph{coupled parabolic systems} attracted the attention of many researchers.
This system class is not only of paramount theoretical appeal, {but it also has a major importance for applications}.
 Typical real-world problems originate from chemical and biochemical engineering, whereby chemical fixed-bed and tubular reactors are important examples (see, e.\:g., \cite{Ra81,Je82}). A first solution to the backstepping design for coupled diffusion-reaction systems with equal diffusion coefficients and constant coefficients was given in \cite{Ba14,Ba15,Ba15a}. In this work it was shown how the approach of \cite{Sm04} for determining the scalar integral kernel, that defines the backstepping transformation, can directly be extended to the matrix case. The latter appears in the control of $n$ coupled PDEs, because the state vector is $n$-dimensional. More challenging is the backstepping control of coupled parabolic systems with distinct diffusion coefficients. First solutions to this problem considered diffusion-reaction systems with constant coefficients (see \cite{Ba15a,LiuB16,Orl17}). For this system class it is possible to simplify the kernel equations by assuming a diagonal structure of the matrix kernel. However, this approach cannot be extended to the spatially-varying coefficients case. Recently, a very general and elegant solution of the backstepping problem for diffusion-convection-reaction systems with spatially-varying coefficients and Dirichlet boundary conditions (BCs) was presented in \cite{Vaz16a}. Thereby, no assumption on the form of integral kernel is imposed. By introducing a local coupling term with a triangular structure in the target system it was shown that the resulting kernel equations are well-posed. The latter result followed from establishing an interesting relation of the kernel equations for the considered system class with the kernel equations related to coupled first-order hyperbolic PDEs treated in \cite{Hu15a,Hu15b}.

A larger class of DPS can be modelled by \emph{partial integro-differential equations (PIDEs)}. They arise directly in the physical modelling or from a singular perturbation of separate subsystems with different time-scales. For this system class, Robin BCs have to be taken into account frequently. This results, e.\:g., from the energy balance along the boundaries in heat transfer problems if the DPS is in direct contact with a surrounding medium. Accordingly, heat isolation at the boundaries requires to introduce Neumann BCs. Moreover, the linearisation of a nonlinear DPS around a steady state, non-homogenous materials or unusually shaped spatial domains lead to spatially-varying coefficients. Hence, it is of interest to extend the backstepping method to this class of coupled parabolic systems.

This paper is concerned with the backstepping control of parabolic PIDEs of Volterra-type with distinct diffusion coefficients and spatially-varying parameters. Thereby, Dirichlet, Neumann or Robin BCs are assumed independently for each component of the state. 
Similar to \cite{Vaz16a} a target system with a local coupling term in triangular form is proposed. This results in a cascaded set of parabolic PDEs, which allow a simple proof of well-posedness and stability. Moreover, the corresponding rate of exponential convergence can directly be  specified in the design so that a simple parametrization of the target system  is possible. The main contribution of the paper is the verification of the well-posedness of the resulting \emph{kernel equations}. For $n$ coupled parabolic PIDEs they are a system of $n^2$ coupled hyperbolic PIDEs with a \emph{Klein-Gordon-type} spatial differential operator, spatially-varying coefficients and rather unusual boundary conditions with integral terms. 
Since this type of kernel equations cannot be traced back to the kernel equations found for general heterodirectional hyperbolic systems in \cite{Hu15a,Hu15b} with the result in \cite{Vaz16a}, a new method for their solution is needed.
In order to provide a systematic solution procedure, the approach of \cite{Sm04,Sm05} is directly extended to the system class in question. The resulting kernel equations for the diagonal elements share the same form as in the scalar case. However, the off-diagonal elements are governed by kernel equations with a different structure.
In a first step, the kernel PIDEs with spatially-varying coefficients are simplified utilizing suitable transformations. In particular, by extending the results in \cite{Sm05}, new kernel PIDEs are obtained, where the coefficients \wrt  the second-order partial derivatives are equal to one. 
Hence, the usual linear change of coordinates mapping the kernel PIDE into its canonical form can directly be applied. The resulting kernel equations allow a very systematic formulation of the corresponding integral equations. This leads to the usual result with the standard triangular spatial domain for the diagonal kernel elements. However, the spatial domain related to the off-diagonal kernel PIDEs is no longer restricted to the first quadrant. It is shown that by introducing suitable artificial BCs, well-posed kernel equations are obtained. The latter property is proved by verifying the uniform convergence of the corresponding successive approximation. This results in a systematic method for the backstepping state feedback stabilization of the coupled PIDEs in question.

As far as coupled systems of parabolic PDEs are concerned, the presented approach provides a direct solution of the $n^2$ second-order hyperbolic kernel PDEs for Dirichlet, Neumann and Robin BCs. In contrast, the recent approach in \cite{Vaz16a} assuming Dirichlet BCs solves a system of $2n^2$ first-order transport equations, which follow from the second-order kernel PDEs. The new approach could also be of interest when considering other types of coupled PDEs. For example, it may be the immediate starting point for an extension to coupled parabolic systems with spatially- and temporally-varying reaction by directly utilizing the results \cite{Sm05,Me09} for scalar parabolic systems.
\vspace{2mm}

The next section formulates the considered problem. Then, the target system is introduced and the corresponding well-posedness and stability is proved in the following section. A systematic procedure for the solution of the kernel equations is presented in Section \ref{sec:solker}. Subsequently, Section \ref{sec:clstab} investigates the stability of the resulting closed-loop system. The proposed backstepping-based state feedback controller design is illustrated for an unstable coupled parabolic system of two PIDEs, where a Dirichlet and Robin BC is imposed at the unactuated boundary and the actuation is of Neumann and Dirichlet type with a coupling at that boundary.
\vfill

\section{Problem formulation}\label{sec:intro}
Consider the following system described by \emph{coupled linear parabolic PIDEs}
\begin{subequations}\label{drs}
 \begin{align}
	\partial_tx(z,t) ={}&%
	 \mathpushleft[0]{\Lambda(z)\partial_z^2x(z,t) + \Phi(z)\partial_zx(z,t) + A(z)x(z,t)}\nonumber\\
	& \mathpushleft[0]{+ A_0(z)x(0,t) + \int_0^zF(z,\zeta)x(\zeta,t)\d\zeta} \label{pdes}\\
	\theta_0[x(t)](0) ={}& 0, && t > 0\label{bc1}\\
	\theta_1[x(t)](1) ={}& u(t),      &&  t > 0 \hspace{4cm}\label{bc2} 	
 \end{align}
\end{subequations}
with \eqref{pdes} defined on $(z,t) \in (0,1) \times \mathbb{R}^+$, the state $x(z,t) \in \mathbb{R}^n$, $n > 1$, and the input $u(t) \in \mathbb{R}^n$. 
The matrix $\Lambda(z) \in \mathbb{R}^{n \times n}$
{given by}
\begin{equation}
\Lambda(z) = \operatorname{diag}(\lambda_1(z),\ldots,\lambda_n(z))
\end{equation}
{contains}
mutually different, positive and spatially-varying \emph{diffusion coefficients} $\lambda_i \in C^{2}[0,1]$, $i = 1,2,\ldots,n$, i.\:e., $\lambda_i(z) \neq \lambda_j(z)$, $z \in [0,1]$, $i \neq j$. The \emph{convection term} is characterized by the diagonal matrix
\begin{equation}\label{phicondef}
\Phi(z) = \operatorname{diag}(\Phi_1(z),\ldots,\Phi_n(z))
\end{equation}
with $\Phi \in (C^1[0,1])^{n \times n}$ and the matrix $A = [A_{ij}]\in (C^1[0,1])^{n \times n}$ describes the \emph{reaction term}. The \emph{local term} in the PIDE \eqref{pdes} is determined by ${A_0} \in (C^1[0,1])^{n \times n}$ and the \emph{non-local term} by the integral kernel $F \in (C^1([0,1]^2))^{n \times n}$. 

The BCs are represented by the formal \emph{matrix differential operators}
\begin{equation}\label{thetadef}
\theta{_i}[h] = B_{{i}}^1\d_zh + B^0_{{i}} h, \quad i = 0,1.
\end{equation}
Therein,
\begin{equation}\label{P1def}
B^1_{0} = \begin{bmatrix}
0_{m} & 0\\
0 & I_{p}
\end{bmatrix} \quad \text{and} \quad
B^0_{{0}} = \begin{bmatrix}
I_{m} & 0\\
0   & Q_0
\end{bmatrix},  
\end{equation}
$m + p = n$, is assumed, which means that the first $m$ BCs at $z = 0$ are of Dirichlet type.  Accordingly, the remaining $p$ BCs represent Neumann/Robin BCs, in which 
\begin{align}\label{eq:Q0}
	Q_0 = \operatorname{diag}(q_1,\ldots,q_{p}) \in \mathbb{R}^{p \times p}
\end{align} is a diagonal matrix. This setup can always be achieved for decoupled BCs at $z = 0$ by a suitable reordering of the state $x$. Moreover, the actuation at $z = 1$ is described by a diagonal matrix $B^1_{1} \in \mathbb{R}^{n \times n}$ and an arbitrary matrix $B^0_{1}\in \mathbb{R}^{n \times n}$. This specifies Dirichlet, Neumann or Robin actuation at $z = 1$ independently for each component $u_i$, $i = 1,2,\ldots,n$, of the input $u$. Thereby, also different types of BCs for each component $x_i(z,t) \in \mathbb{R}$, $i = 1,2,\ldots,n$, of the state $x$ at $z = 1$ as well as possible couplings of these BCs are included.
\begin{rem}
The type of convection and BCs specified by \eqref{phicondef} and \eqref{thetadef} appear frequently in applications. Well-known examples are chemical fixed-bed and tubular reactors (see, e.\:g., \cite{Ra81,Je82}). \hfill $\triangleleft$	
\end{rem}

The initial conditions (ICs) of the system are $x(z,0) = x_0(z)$ with $x_0 \in (L_2(0,1))^n$.

The system \eqref{drs} can be simplified by introducing the boundedly invertible \emph{Hopf-Cole-type state transformation}
\begin{equation}\label{hctraf}
\check{x}(z,t) = \exp\Big(\tfrac{1}{2}\int_0^z\Lambda^{-1}(\zeta)\Phi(\zeta)\d \zeta\Big)x(z,t).
\end{equation}
With a direct calculation it is easy to verify that \eqref{hctraf} maps \eqref{pdes} into the  coupled PIDEs
\begin{align}
 \partial_t\check{x}(z,t) &= \Lambda(z)\partial_z^2\check{x}(z,t) + \check{A}(z)\check{x}(z,t)\nonumber\\
 & \quad  + \check{A}_0(z)\check{x}(0,t) + \int_0^z\check{F}(z,\zeta)\check{x}(\zeta,t)\d\zeta\label{dcrpde}
\end{align}
on $(z,t) \in (0,1) \times \mathbb{R}^+$. Thereby, the resulting system shares the same type of BCs as \eqref{drs}. Consequently, the convection term in \eqref{drs} is omitted in the sequel, i.\:e., $\Phi(z) \equiv 0$ is assumed.

Consider the \emph{state feedback controller}
\begin{equation}\label{sfeed}
u(t) = \mathcal{K}[x(t)]
\end{equation}
with the formal \emph{feedback operator} $\mathcal{K}$ mapping the state $x(z,t) \in \mathbb{R}^n$ to the input $u(t) \in \mathbb{R}^n$. The problem considered in this paper is the backstepping design of \eqref{sfeed} such that the resulting closed-loop system is exponentially stable with a prescribed rate of convergence.

\section{State feedback controller design}\label{sec:kerder}
\subsection{Selection of the Target System}
In what follows, the \emph{backstepping approach} (see, e.\:g., \cite{Kr08}) is used to determine the state feedback controller \eqref{sfeed}. To this end, the \emph{backstepping coordinates}
\begin{equation}\label{btrafo}
\tilde{x}(z,t) = x(z,t) - \int_0^zK(z,\zeta)x(\zeta,t)\d \zeta
\end{equation}
with the \emph{integral kernel} $K(z,\zeta) \in \mathbb{R}^{n \times n}$ are introduced for \eqref{drs}. The determination of the feedback gains in \eqref{sfeed} {is} significantly simplified in these coordinates. As a first step in the design, an exponentially stable target system has to be found so that the change of coordinates \eqref{btrafo} exists. By following the lines of \cite{Sm04} for the scalar case and the recent results in \cite{Hu15a,Vaz16a}, the following \emph{target system}
\begin{subequations}\label{drst}
	\begin{align} \label{tpdes}
	\partial_t\tilde{x}(z,t) &= \Lambda(z)\partial_z^2\tilde{x}(z,t)  - \mu_c\tilde{x}(z,t)\\ \nonumber
	&\quad  - \widetilde{A}_0(z)(E_1E_1^{\top}\partial_z\tilde{x}(0,t) + E_2E_2^{\top}\tilde{x}(0,t)) \\
	\theta_0[\tilde{x}(t)](0) &= 0, &&\hspace{-5cm}  t > 0\label{tbc1}\\
	\tilde{\theta}_{1}[\tilde{x}(t)](1) &= 0,  &&\hspace{-5cm} t > 0\label{tbc2}
\end{align}
\end{subequations}
with $\mu_c \in \mathbb{R}$ is proposed. Therein, \eqref{tpdes} is defined on $(z,t) \in (0,1) \times \mathbb{R}^+$ and 
$\widetilde{A}_{0}(z) = [\widetilde{A}_{0,ij}(z)] \in \mathbb{R}^{n \times n}$ 
is a matrix containing only $n(n-1)/2$ non-zero elements, \ie
\begin{align} \label{eq:A0_til}
\widetilde{A}_{0,ij}(z) = 
\begin{cases}
f_{ij}(z), & \lambda_i < \lambda_j  \\
0, & \text{else}
\end{cases}
\end{align}
with $f_{ij}(z)$ determined by the kernel $K(z,\zeta)$.
Here, $\lambda_i < \lambda_j$ is the shorthand notation for $\lambda_i(z) < \lambda_j(z)$, $z \in [0,1]$.
%
The non-zero elements $\widetilde{A}_{0,ij}(z)$, $\lambda_i < \lambda_j$, are determined by the solution of the kernel equations (cf. Section \ref{sec:cform}). Furthermore, the matrices
\begin{equation}
E_1 = \begin{bmatrix}
I_{m}\\
0
\end{bmatrix}  \in \mathbb{R}^{n \times m} \quad \text{and} \quad
E_2 = \begin{bmatrix}
0\\
I_{p}
\end{bmatrix} \in \mathbb{R}^{n \times p}
\end{equation}
were used in \eqref{tpdes}. The BCs of the target system {at $z = 1$} are characterized by
$\tilde{\theta}_{1}[h] = \widetilde{B}_{1}^1\d_zh + \widetilde{B}_{1}^0 h$ with
\begin{equation}
 \widetilde{B}_{1}^1 = \operatorname{diag}({\tilde{b}}^{1}_1,\ldots,{\tilde{b}}^{1}_n) \quad \text{and} \quad \widetilde{B}_{1}^0 = \operatorname{diag}({\tilde{b}}^{0}_1,\ldots,{\tilde{b}}^{0}_n)
\end{equation}
being diagonal matrices satisfying $|\tilde{b}_i^{1}| + |\tilde{b}_i^{0}| > 0$, $i = 1,2,\ldots,n$. This means that decoupled BCs at $z = 1$ are assigned in the target system. 
Thereby, Dirichlet BCs are assigned for $z=1$, if the plant has this type of BC for the considered element of $x$. 
On the contrary,
Neumann or Robin BCs can be interchanged. 
As there are $p$ Neumann/Robin BCs in \eqref{bc1},
the local coupling term in \eqref{tpdes} has to depend on $\tilde{x}(0,t)$. In contrast, for Dirichlet BCs it was shown in \cite{Vaz16a} that the corresponding term requires a dependency on $\partial_z\tilde{x}(0,t)$. These couplings have to be introduced, because a complete decoupling of the PDEs in \eqref{tpdes} leads to an overdetermined and thus unsolvable set of kernel equations.

For the stability analysis of the target system \eqref{drst} it is convenient to introduce the new state
\begin{equation}
 \tilde{x}^*(z,t) = P\tilde{x}(z,t).
\end{equation}
Therein, $P \in \mathbb{R}^{n \times n}$ is a \emph{permutation matrix} with $P^{-1} = P^{\top}$ such that  the relation $\lambda^*_1(z) > \ldots > \lambda^*_n(z)$  holds for the diagonal elements in $\Lambda^*(z) = \operatorname{diag}(\lambda^*_1(z),\ldots,\lambda^*_n(z)) = P\Lambda(z)P^{\top}$. Then, the target system \eqref{drst} takes the form%
\begin{subequations}\label{drst2}%
\begin{align}\label{tpdes2}%
	\partial_t\tilde{x}^*(z,t) &= \Lambda^*(z)\partial_z^2\tilde{x}^*(z,t)  - \mu_c\tilde{x}^*(z,t)\\ \nonumber
	&\quad  - \widetilde{A}^*_0(z)(R_1\partial_z\tilde{x}^*(0,t) + R_2\tilde{x}^*(0,t)) \\
	\theta_0[P^{\top}\tilde{x}^*(t)](0) &= 0 \label{tbc12}\\
	\tilde{\theta}_{1}[P^{\top}\tilde{x}^*(t)](1) &= 0 \label{tbc22}
	\end{align}
\end{subequations}
with $R_1 = PE_1E_1^{\top}P^{\top}$ and $R_2 = PE_2E_2^{\top}P^{\top}$. Due to the previous definition \eqref{eq:A0_til} of $\widetilde{A}_{0}(z)$ and the introduced permutation of the state, the matrix $\widetilde{A}^*_{0}(z)$ is strictly lower triangular, i.\:e.,
\begin{equation}\label{a01def}
\widetilde{A}^*_{0}(z) = P\widetilde{A}_{0}(z)P^{\top}
= 
\begin{bmatrix} 0     & \ldots &  \ldots           & 0\\
\widetilde{A}^*_{0,21}(z) & \ddots       & \ddots             & \vdots\\
\vdots    & \ddots &  \ddots     &  \vdots \\
\widetilde{A}^*_{0,n1}(z) & \ldots & \widetilde{A}^*_{0,n\,n-1}(z) & 0
\end{bmatrix}.
\end{equation}
This implies that \eqref{drst2} has a cascade structure, which significantly simplifies the corresponding stability analysis. The next theorem shows that the target system \eqref{drst2} and thus \eqref{drst} is exponentially stable with a stability margin assignable by $\mu_c$.
\begin{thm}[Stability of the target system]\label{thm:tstab}
Assume that $\mu_c > \mu_{\text{max}}$, in which $\mu_{\text{max}}$ is the largest eigenvalue of \eqref{drst2} for $\mu_c = 0$ and $\widetilde{A}^*_0(z) \equiv 0$. Then, the target system \eqref{drst2} is exponentially stable in the $L_2$-norm $\|h\| = (\int_0^1\|h(z)\|^2_{\mathbb{C}^n}\d z)^{1/2}$, i.\:e.,
\begin{equation}\label{expstab}%
\|\tilde{x}(t)\| \leq \widetilde{M}\text{e}^{(\mu_{\text{max}}-\mu_c) t}\|\tilde{x}(0)\|,\quad t \geq 0
\end{equation}%
for all $\tilde{x}(0) \in (H^2(0,1))^n$ compatible with the BCs \eqref{tbc12}, \eqref{tbc22} and an $\widetilde{M} \geq 1$. 
\end{thm}	
The proof of this result can be found in the appendix. 

\subsection{Derivation of the Kernel Equations}
\label{sec:keq}
The equations to be solved for determining $K(z,\zeta)$ in \eqref{btrafo} result from requiring that \eqref{btrafo} and a suitable feedback \eqref{sfeed} map \eqref{drs} into the target system \eqref{drst}. Differentiating \eqref{btrafo} \wrt time, inserting \eqref{pdes}, utilizing \eqref{tpdes} and interchanging the order of integration in the double integral results in
\begin{align}\notag
    & \partial_t\tilde{x}(z,t) = \Lambda(z)\partial_z^2\tilde{x}(z,t)
		- \mu_c\tilde{x}(z,t)\\\notag
	&  - \widetilde{A}_0(z)(E_1E_1^{\top}\partial_z\tilde{x}(0,t) + E_2E_2^{\top}\tilde{x}(0,t))\\\notag
	& + \Lambda(z)\partial_z^2\int_0^zK(z,\zeta)x(\zeta,t)\d \zeta\\\notag
	& + \Big(\widetilde{A}_0(z)E_2E_2^{\top} + A_0(z) - \int_0^zK(z,\zeta){A_0}(\zeta)\d\zeta\Big) x(0,t)\\\notag
	&  +\widetilde{A}_0(z)E_1E_1^{\top}\partial_zx(0,t) + (A(z) + \mu_cI)x(z,t) \\\notag
	&
	 -\int_0^zK(z,\zeta)\Lambda(\zeta)\partial_{\zeta}^2x(\zeta,t)\d \zeta  -\int_0^z\Big(K(z,\zeta)(A(\zeta) + \mu_cI)\\
	& + \int_{\zeta}^{z}K(z,\ozeta)F(\ozeta,\zeta)\d\ozeta \Big)x(\zeta,t)\d\zeta,
\end{align}
in which $K(0,0) = 0$ was assumed. It is convenient to define the integral operators
\begin{multline}
	\mathcal{B}[K](z,\zeta) = K(z,\zeta)(A(\zeta) + \mu_cI)\\
	- F(z,\zeta) + \int_{\zeta}^{z}K(z,\ozeta)F(\ozeta,\zeta)\d\ozeta
\end{multline}
and
\begin{equation}
	\mathcal{C}[K](z) ={ A_0}(z) - \int_0^zK(z,\zeta){A_0}(\zeta)\d\zeta.
\end{equation}
Then, integration by parts and using the Leibniz differentiation rule yield after simple calculations 
\begin{align}\notag
	 \partial_t \tilde{x}(z,t)  ={}& \Lambda(z)\partial_z^2\tilde{x}(z,t) - \mu_c\tilde{x}(z,t)\\ \notag
	& - \widetilde{A}_0(z)\big(E_1E_1^{\top}\partial_z\tilde{x}(0,t) + E_2E_2^{\top}\tilde{x}(0,t)\big)\\\notag
	& + \int_0^z\Big(\Lambda(z) \dz^2 K(z,\zeta) - {\dzeta^2}(K(z,\zeta)\Lambda(\zeta)) \\\notag
	&  - \mathcal{B}[K](z,\zeta)\Big)x(\zeta,t)\d \zeta\\\notag
	& + M_1(z)x(0,t) + M_2(z)\partial_z x(0,t)\\\notag
	& + \Big(\Lambda(z) K{'}(z,z) + \Lambda(z) \dz K(z,z)\\\notag
	& + \dzeta K(z,z)\Lambda(z) + K(z,z) \Lambda{'}(z) \\\notag
	& +  A(z) + \mu_cI\Big)x(z,t) \\
	& + \big(\Lambda(z)K(z,z) - K(z,z)\Lambda(z)\big)\partial_zx(z,t) \label{newcoordsys}
\end{align}%
with 
{$(\,\cdot\,){'} = \frac{\d}{\d z}(\,\cdot\,)$}.
Therein, $M_1(z) = \widetilde{A}_0(z)E_2E_2^{\top} + \mathcal{C}[K](z) - \dzeta K(z,0)\Lambda(0)  - K(z,0)\Lambda'(0)$ and $M_2(z) = \widetilde{A}_0(z)E_1E_1^{\top} + K(z,0)\Lambda(0)$ were utilized.
In order to derive the kernel BCs use $E_1E_1^{\top} + E_2E_2^{\top} = I_n$ and ${x^i} = E_i^{\top}x$, $i = 1,2$, as well as the condition
\begin{align}\notag
& M_1(z)x(0,t) + M_2(z)\partial_zx(0,t)\\ \notag
& = M_1(z)E_1x{^1}(0,t) + M_1(z)E_2x{^2}(0,t)\\
& \quad + M_2(z)E_1\partial_z x{^1}(0,t) + M_2(z)E_2\partial_z x{^2}(0,t) = 0 \label{z0bed}
\end{align}
implied by {the requirement that \eqref{newcoordsys} coincides with \eqref{tpdes}}. 
Observe that \eqref{bc1} gives $x{^1}(0,t) = 0$ and $\partial_z x{^2}(0,t) = -Q_0 x{^2}(0,t)$
in view of \eqref{thetadef} and \eqref{P1def}. After inserting this in \eqref{z0bed} and comparing the BC \eqref{bc1} with \eqref{tbc1}, the \emph{kernel equations}
\begin{subequations}\label{keq}
	\begin{align}
	\tikzmark[0pt]{first}{}
	 \Lambda(z) \dz^2 K(z,\zeta) - \dzeta^2 (K(z,\zeta)\Lambda(\zeta))  
	 &= \mathcal{B}[K](z,\zeta)\label{kpde}\\
	 K(z,0)\Lambda(0)E_1 &= -\widetilde{A}_0(z)E_1\label{dtbc}\\[4pt]
	\tikzmark[0pt]{line}{}\strut\label{kBC1}
	\\\notag
	&= \!\big(\widetilde{A}_0(z) \!+\! \mathcal{C}[K](z)\big)E_2\\[4pt]
	\tikzmark[0pt]{line2}{}\strut \notag
	\\\label{kBC2}
	+ K(z,z)\Lambda'(z)  &= -(A(z) \!+\! \mu_cI)\\[4pt]
	K(z,z)\Lambda(z) - \Lambda(z)K(z,z) &= 0\label{kBC3}\\  
	K(0,0) &= 0\label{kIC}
	\tikz[overlay, remember picture]{
	 \node at (first|-line)[anchor=base west, inner sep=0pt]{$\dzeta K(z,0)\Lambda(0)E_2 + K(z,0)(\Lambda'(0)E_2 + \Lambda(0)E_2Q_0)$};
	 \node at (first|-line2)[anchor=base west, inner sep=0pt]
	 {$\Lambda(z) K'(z,z) \!+\! \Lambda(z) \dz K(z,z) \!+\! \dzeta K(z,z)\Lambda(z)$};
	}
	\end{align}
\end{subequations}
are obtained, in which \eqref{kpde} is defined on $0 < \zeta < z < 1$.
\begin{rem}
	If all BCs at $z = 0$ are of the same type, then the substitution $E_1 \to I_n$ and $E_2 \to 0$ for Dirichlet BCs and  $E_1 \to 0$ and $E_2 \to  I_n$ for Neumann/Robin BCs in \eqref{keq} yields the related kernel equations. With this, their solution also follows from the subsequent results. \hfill $\triangleleft$	
\end{rem}

If the kernel $K(z,\zeta)$ is known, then the feedback controller \eqref{sfeed} follows from imposing the BC \eqref{tbc2}. For this, consider \eqref{bc2} with \eqref{thetadef} in the form
\begin{equation}\label{upara}
 u(t) = \theta_1[x(t)](1) = B_{1}^1\partial_zx(1,t) + B_1^0 x(1,t).
\end{equation}
If an element in the coefficient matrix $\widetilde{B}_{1}^1$ of  $\tilde{\theta}_1[\tilde{x}(t)](1)$ (see \eqref{tbc2}) is zero, \ie $\tilde{b}^1_i = 0$,
then $\tilde{x}_i(1,t) = 0$ and \eqref{btrafo} yields {for the $i$-th component of the state vector}
\begin{equation}\label{rel1}
 x_i(1,t) = \int_0^1e_i^{\top}K(1,z)x(z,t)\d z.
\end{equation} 
Therein, $e_i \in \mathbb{R}^n$ denotes the $i$-th unit vector. 
In the case $\tilde{b}^1_i \neq 0$, one has $\dz\tilde{x}_i(1,t) = -(\tilde{b}_i^0/\tilde{b}_i^1) \tilde{x}_i(1,t)$ and obtains
\begin{multline}\label{rel2}
 \partial_zx_i(1,t) = (e_i^{\top}K(1,1) - \frac{{\tilde{b}}^{0}_i}{{\tilde{b}}^{1}_i}e_i^{\top})x(1,t)\\
  + \int_0^1e_i^{\top}({\dz} K(1,z) + \frac{{\tilde{b}}^{0}_i}{\tilde{b}^{1}_i}K(1,z))x(z,t)\d z
\end{multline}
after differentiating \eqref{btrafo} \wrt $z$.
With this, the feedback operator $\mathcal{K}$ in \eqref{sfeed} follows from inserting \eqref{rel1} and \eqref{rel2} in \eqref{upara}. Thereby, all spatial derivatives of the state possibly appearing in \eqref{upara} are removed by the substitution \eqref{rel2}.

In the next section the solvability of \eqref{keq} is investigated. To this end, \eqref{keq} is converted into integral equations which can be solved by a fixpoint iteration. By considering a sufficiently large but finite number of iterations this leads to the \emph{method of successive approximations}, which results in  a systematic approach for determining the kernel $K(z,\zeta)$. 

\section{Solution of the kernel equations}\label{sec:solker}
\subsection{Component Form of the Kernel Equations}\label{sec:cform}
For converting \eqref{keq}  into integral equations the boundary value problems (BVPs) for the elements $K_{ij}(z,\zeta)$ of the kernel $K(z,\zeta) = [K_{ij}(z,\zeta)]$, $i,j = 1,2,\ldots,n$, satisfying \eqref{keq} are derived. Thereby, $\mathcal{B} = [\mathcal{B}_{ij}]$ and $\mathcal{C} = [\mathcal{C}_{ij}]$ are utilized. This yields the \emph{kernel PIDEs}
\begin{equation}\label{kpdeineqj}
\lambda_i(z)\partial_z^2K_{ij}(z,\zeta) \!-\! \partial^2_{\zeta}(\lambda_j(\zeta)K_{ij}(z,\zeta))
 = \mathcal{B}_{ij}[K](z,\zeta).
\end{equation}
Evaluating \eqref{dtbc} for $\lambda_i \geq \lambda_j$ and $j \leq m$ leads to the BC
\begin{equation}\label{dmatbcdir}
K_{ij}(z,0) = 0
\end{equation}
in view of \eqref{eq:A0_til}.
Similarly, 
for $\lambda_i \geq \lambda_j$ and $j > m$, \eqref{kBC1} with \eqref{eq:Q0} results in the BC
\begin{equation}\label{dmatbc}
\lambda_j(0)\partial_{\zeta}K_{ij}(z,0) + (\lambda'_j(0) + q_j\lambda_j(0))K_{ij}(z,0)
 = \mathcal{C}_{ij}[K](z),
\end{equation}
whereas the relations
\begin{subequations}\label{eq:Aij_tilde}
\begin{empheq}[left={\!\widetilde{A}_{0,ij}(z) =\empheqlbrace}]%
{align}%
 &-\lambda_j(0)K_{ij}(z,0), &&\hspace{-2mm} j \leq m \\[10pt] \label{Aijtildet}%
  &\begin{aligned}[c]
 & \lambda_j(0)\partial_{\zeta}K_{ij}(z,0)  - \mathcal{C}_{ij}[K](z)%
 \\&+ \big(\lambda'_j(0) + q_j\lambda_j(0)\big)K_{ij}(z,0), \hspace{-1mm}  
  \end{aligned}  &&\hspace{-2mm}  j > m
\end{empheq}
\end{subequations}%
%
obtained for $\lambda_i < \lambda_j$ determine the elements of $\widetilde{A}_0(z)$ in \eqref{eq:A0_til}. 
Next, considering \eqref{kBC3} for $i = j$ results in
\begin{equation}
(\lambda_i(z) - \lambda_i(z))K_{ii}(z,z) = 0,
\end{equation}
which shows that no BC has to be fulfilled for $K_{ii}(z,z)$ {to satisfy \eqref{kBC3}}. In contrast, the case $i \neq j$ gives
\begin{equation}\label{zzbed}
(\lambda_j(z) - \lambda_i(z))K_{ij}(z,z) = 0
\end{equation}
so that for unequal diffusion coefficients $\lambda_i(z)$ and $\lambda_j(z)$ the BC
\begin{equation}\label{kijzero}
K_{ij}(z,z) = 0
\end{equation}
follows. 
After differentiating \eqref{kijzero} one obtains 
\begin{equation}\label{dzkij}
K_{ij}'(z,z) = 0
\end{equation}
and thus
\begin{equation}\label{pdiffcom}
\partial_{\zeta}K_{ij}(z,z) = -\partial_zK_{ij}(z,z).
\end{equation}
Finally, the remaining BC \eqref{kBC2} for $i = j$ leads to the ODE
\begin{equation}
2\lambda_i(z)K_{ii}'(z,z) + \lambda'_i(z)K_{ii}(z,z) = -(A_{ii}(z) + \mu_c).
\end{equation}
With the IC \eqref{kIC}{,} the corresponding solution is
\begin{equation}
K_{ii}(z,z) = -\frac{1}{\sqrt{\lambda_i(z)}}\int_0^z\frac{A_{ii}(\zeta) + \mu_c}{2\sqrt{\lambda_i(\zeta)}}\d \zeta.
\end{equation}
If $i \neq j$, then \eqref{kBC2} gives
\begin{equation}
(\lambda_i(z) - \lambda_j(z))\partial_zK_{ij}(z,z) = - A_{ij}(z)
\end{equation}
in view of \eqref{dzkij} and \eqref{pdiffcom}. Hence, the BC
\begin{equation}
\partial_zK_{ij}(z,z) = \frac{A_{ij}(z)}{\lambda_j(z) - \lambda_i(z)}
\end{equation}
follows. By collecting the preceding results{,} the following kernel equations for the elements of $K(z,\zeta)$ can be deduced. Thereby, the BC in $[\,\cdot\,]_{\ast}$ has to be considered only for those elements of the corresponding expression, for which the condition $\ast$ is fulfilled. With this, the \emph{component form}
\begin{subequations}\label{keqcompdiag}
	\begin{flalign}
	& \mathpushleft[0]{\underline{i = j:}}\nonumber\\[0.1cm]
	& \mathpushleft[0]{\lambda_i(z)\partial_z^2K_{ii}(z,\zeta) - \partial^2_{\zeta}(\lambda_i(\zeta)K_{ii}(z,\zeta)) = \mathcal{B}_{ii}[K](z,\zeta)} 
	&& \phantom{(z,\zeta) - \partial^2_{\zeta}(\lambda_i(\zeta)K_{ii}(z,\zeta)) = \mathcal{B}_{ii}[K](z,\zeta)} 
	\label{keqcompdiagpde}\\
	&&[K_{ii}(z,0) &= 0]_{i \leq m}\\
	&\mathpushleft[0]{[\lambda_i(0)\partial_{\zeta}K_{ii}(z,0) \! + \!(\lambda'_i(0) \!+\! q_i\lambda_i(0))K_{ii}(z,0) }\nonumber\\
	&&& = \mathcal{C}_{ii}[K](z)]_{i > m}\\
	&& K_{ii}(z,z) &= -\frac{1}{\sqrt{\lambda_i(z)}}\int_0^z\frac{A_{ii}(\zeta) + \mu_c}{2\sqrt{\lambda_i(\zeta)}}\d \zeta\label{zzklassbed}
	\end{flalign}
\end{subequations}
and
\begin{subequations}\label{keqcomp1}
	\begin{flalign}
	&\mathpushleft[0]{\underline{i \neq j:}}\nonumber\\[0.1cm]
	&\mathpushleft[0]{\lambda_i(z)\partial_z^2K_{ij}(z,\zeta) - \partial^2_{\zeta}(\lambda_j(\zeta)K_{ij}(z,\zeta)) = \mathcal{B}_{ij}[K](z,\zeta)}
\label{keqcomp1pde} \\
	&&[K_{ij}(z,0) &= 0]_{\lambda_i > \lambda_j,\, j \leq m}\\
	&\mathpushleft[0]{[\lambda_j(0)\partial_{\zeta}K_{ij}(z,0) + (\lambda'_j(0) + q_j\lambda_j(0))K_{ij}(z,0)} \nonumber\\
	&&&= \mathcal{C}_{ij}[K](z)]_{\lambda_i > \lambda_j,\, j > m}\label{bbc}\\
	&& K_{ij}(z,z) &= 0\label{zzbedkeq}\\
	&& \partial_zK_{ij}(z,z) &= \frac{A_{ij}(z)}{\lambda_j(z) - \lambda_i(z)}\label{zzpzbed} \hspace{3.2cm}
	\end{flalign}
\end{subequations}
of the kernel equations \eqref{keq} is obtained, whereby \eqref{keqcompdiagpde} and \eqref{keqcomp1pde} are defined on $0 < \zeta < z < 1$.
\begin{rem}
It is interesting to note that the kernel equations for the diagonal elements \eqref{keqcompdiag} consist of the usual but coupled kernel equations found in the scalar case (see \cite{Sm04,Sm05}). In contrast, the BVP \eqref{keqcomp1} for all other elements are different \wrt the kernel PIDE and the corresponding BCs, which needs a new approach for their solution. This is the topic of the next sections. \hfill $\triangleleft$
\end{rem}

\begin{rem}
In the case of \emph{equal diffusion coefficients}, i.\:e., $\lambda_{\text{max}} \geq  \lambda_1(z) = \lambda_2(z) = \ldots = \lambda_{n}(z) \geq \lambda_{\text{min}} > 0$, $z \in [0,1]$, the kernel equations \eqref{keqcomp1} have the same form as \eqref{keqcompdiag} and thus are easier to solve. In particular, \eqref{zzbed} leads to no condition for $K_{ij}(z,z)$, i.\:e., \eqref{zzbedkeq} does not appear. Then, \eqref{kBC2} yields an ODE for $K_{ij}(z,z)$ replacing \eqref{zzpzbed} by a BC of the form \eqref{zzklassbed}. As a consequence, \eqref{dmatbcdir} and \eqref{dmatbc} can be imposed for all $i \neq j$ as kernel BCs so that $\widetilde{A}_0(z) \equiv 0$ holds in \eqref{tpdes}. This result was first shown in \cite{Ba14} for a diffusion-reaction system with constant coefficients and Neumann BCs. 
\hfill $\triangleleft$
\end{rem}

\subsection{Transformation of the Kernel Equations}\label{sec:ctraf}
In what follows, the transformation approach in \cite{Sm05} for a single kernel PIDE is extended to the coupled system of kernel PIDEs \eqref{keqcompdiagpde} and \eqref{keqcomp1pde}. Then, the well-posedness of the BVPs related to the resulting simpler kernel PIDEs is much easier to prove. 

For the ease of presentation, {matrix elements and coordinates} with indices $i$ and $j$ may be represented without the index but bold face, e.\:g., $\pK = K_{ij} $, holds in the sequel. Moreover, index notation for derivatives is used, \ie $f_z = \partial_z f$.
 
\subsubsection{Transformation of the second-order partial derivatives}
In order to eliminate the dependency on $z$ and $\zeta$ of the coefficients \wrt the second-order partial derivatives in \eqref{keqcompdiagpde} and \eqref{keqcomp1pde}, introduce the change of coordinates
\begin{align}\label{cchange}
\bm{\rho} = \rho{_i}(z) = \phi_i(z) \quad\quad \text{and} \quad\quad \bm{\sigma} =\sigma{_j}(\zeta) =  \phi_j(\zeta).
\end{align}
For notational clarity, $(\prho,\psigma)$ denotes a point in the new coordinate system, whereas $\rho_i(z)=\phi_i(z)$ and $\sigma_j(\zeta)=\phi_j(\zeta)$ are the respective transformations.

Further, define
\begin{align}\label{pbardef}
\bar{\bm{K}}(\rho{_i}(z),\sigma{_j}(\zeta)) = \lambda_j(\zeta)\bm{K}(z,\zeta).
\end{align}
Then, differentiating \eqref{pbardef} twice \wrt $z$ and $\zeta$ as well as inserting the result in \eqref{keqcompdiagpde} and \eqref{keqcomp1pde} yields
\begin{align}\label{prenewtrans}
&\lambda_i(z)(\phi_i'(z))^2\bar{\bm{K}}_{\bm{\rho}\bm{\rho}}(\bm{\rho},\bm{\sigma}) - \lambda_j(\zeta)(\phi_j'(\zeta))^2\bar{\bm{K}}_{\bm{\sigma}\bm{\sigma}}(\bm{\rho},\bm{\sigma})\nonumber\\
& +\lambda_i(z)\phi_i''(z)\bar{\bm{K}}_{\bm{\rho}}(\bm{\rho},\bm{\sigma}) - \lambda_j(\zeta)\phi_j''(\zeta)\bar{\bm{K}}_{\bm{\sigma}}(\bm{\rho},\bm{\sigma})\nonumber\\
& = \lambda_j(\zeta)\bm{\mathcal{B}}[K](z,\zeta)
\end{align}
after a simple calculation. Hence, one obtains the \emph{standard form}
\begin{multline}\label{prenewtrans2}
\bar{\bm{K}}_{\bm{\rho}\bm{\rho}}(\bm{\rho},\bm{\sigma}) -   \bar{\bm{K}}_{\bm{\sigma}\bm{\sigma}}(\bm{\rho},\bm{\sigma}) +\lambda_i(z)\phi_i''(z)\bar{\bm{K}}_{\bm{\rho}}(\bm{\rho},\bm{\sigma})\\ - \lambda_j(\zeta)\phi_j''(\zeta)\bar{\bm{K}}_{\bm{\sigma}}(\bm{\rho},\bm{\sigma})
 = \lambda_j(\zeta)\bm{\mathcal{B}}[K](z,\zeta)
\end{multline}
of the kernel PIDE, if
\begin{equation}\label{phibed}
	\lambda_i(z)(\phi_i'(z))^2 = 1,\quad z \in (0,1), \quad i = 1,2,\ldots,n
\end{equation}
holds. It is straightforward to verify that a solution of \eqref{phibed} is
\begin{align}\label{phidef}
\phi_i(z) = \int_0^{z}\frac{\d\zeta}{\sqrt{\lambda_i(\zeta)}},\quad z \in [0,1],\quad i = 1,2,\ldots,n.
\end{align}
Note that \eqref{phidef} also determines the function $\phi_j(\zeta)$ by appropriate substitution{s}.
The function $\phi_i$ is invertible, because it is strictly monotonically increasing. Hence, the inverse change of coordinates
\begin{subequations}\label{eq:inv_cchange}
	\begin{align}
		z &= z{_i}(\bm{\rho}) = \phi_i^{-1}(\prho)\\
		\zeta &= \zeta{_j}(\bm{\sigma}) = \phi_j^{-1}(\psigma)
	\end{align}
\end{subequations}
related to \eqref{cchange} exists. 
Thereby, $(z,\zeta)$ denotes a point in the original coordinate system, whereas $z_i(\prho) = \phi_i^{-1}(\prho)$ and $\zeta_j(\psigma) = \phi_j^{-1}(\psigma)$ are the corresponding transformations.
 
\subsubsection{Elimination of the first-order partial derivatives} 
In view of \eqref{phidef} the result
\begin{equation}
 \phi_i''(z) = -\frac{\lambda_i'(z)}{2\lambda_i(z)\sqrt{\lambda_i(z)}}
\end{equation}
can be deduced by differentiation.
With this, the kernel PIDEs \eqref{prenewtrans2} can be rewritten as
\begin{align}\label{prenewtrans3}
	&\bar{\bm{K}}_{\bm{\rho}\bm{\rho}}(\bm{\rho},\bm{\sigma}) 
	- \bar{\bm{K}}_{\bm{\sigma}\bm{\sigma}}(\bm{\rho},\bm{\sigma}) \nonumber \\
	& - \kappa_i(z{_i}(\bm{\rho}))
		\bar{\bm{K}}_{\bm{\rho}}(\bm{\rho},\bm{\sigma}) 
	+ \kappa_j(\zeta{_j}(\bm{\sigma}))
		\bar{\bm{K}}_{\bm{\sigma}}(\bm{\rho},\bm{\sigma})\nonumber\\
	& = \lambda_j(\zeta_j(\psigma))\bm{\mathcal{B}}[K](z{_i}(\prho),\zeta{_j}(\psigma)),
\end{align}
in which
\begin{equation}\label{coefffunc}
\kappa_i(z) = \frac{\lambda_i'(z)}{2\sqrt{\lambda_i(z)}} 
\end{equation}
and $\kappa_j(\zeta)$ follows from the corresponding substitutions.
In view of \eqref{keqcompdiagpde},  \eqref{keqcomp1pde} and \eqref{prenewtrans3} the {nonlinear} transformation \eqref{cchange} introduces first-order partial derivatives in the kernel PIDEs. They, however, can readily be eliminated, because the corresponding coefficient functions $\kappa_i(z{_i}(\bm{\rho}))$ and $\kappa_j(\zeta{_j}(\bm{\sigma}))$ depend only on a single variable. Towards this end, the change of coordinates
\begin{equation}\label{secondcc}
\bar{\bm{K}}(\bm{\rho},\bm{\sigma}) = \psi_i(\bm{\rho})\psi_j(\bm{\sigma})\widetilde{\pmb{K}}(\bm{\rho},\bm{\sigma})
\end{equation}
is introduced. From this and \eqref{pbardef} the relation
\begin{align} \label{eq:K_til_orig}
	\pK(z_i(\prho),\zeta_j(\psigma)) = \frac{\psi_i(\prho)\psi_j(\psigma)}{\lambda_j(\zeta_j(\psigma))}\widetilde{\pK}(\prho,\psigma)
\end{align}
can be deduced.
A straightforward but lengthy calculation shows that \eqref{prenewtrans3} takes the form
\begin{equation}\label{prenewtrans4}
\widetilde{\pmb{K}}_{\bm{\rho}\bm{\rho}}(\bm{\rho},\bm{\sigma}) - \widetilde{\pmb{K}}_{\bm{\sigma}\bm{\sigma}}(\bm{\rho},\bm{\sigma})
 =  \widetilde{\pmb{\mathcal{B}}}[\widetilde{K}](\bm{\rho},\bm{\sigma})
\end{equation}
with
\begin{align}\label{eq:B} 
&\widetilde{\pmb{\mathcal{B}}}[\widetilde{K}](\bm{\rho},\bm{\sigma}) =
\Big(
	\pa\big(z{_i}(\prho),\zeta{_j}(\psigma) \big) + \mu_c 
\Big)\widetilde{\pmb{K}}(\prho,\psigma) \\ \notag
&-\underbrace{
	\tfrac{\lambda_j(\zeta{_j}(\psigma))}{\psi_i(\prho)\psi_j(\psigma)}
	\bm{F}(z{_i}(\prho),\zeta{_j}(\psigma))}
	_{\bm{c}^{1}(\prho,\psigma)}   \\ \notag
&
+ \sum\limits_{k=1}^{n}  
\widetilde{K}_{ik}(\prho,\sigma_k(\zeta_j(\psigma)))
\underbrace{
	A_{kj}(\zeta{_j}(\psigma))
	\tfrac{\lambda_j(\zeta_j(\psigma)) \psi_k(\psigma)}{\lambda_k(\zeta_j(\psigma))\psi_j(\psigma)} 
	}_
	{c^2_{kj}(\psigma)}
\\ \notag
& + \int_{\zeta{_j}(\psigma)}^{z_i(\prho)}
\widetilde{K}_{ik}(\prho,\sigma_k(\ozeta))
\underbrace{
		F_{kj}(\ozeta,\zeta{_j}(\psigma))
		\tfrac{\lambda_j(\zeta{_j}(\psigma))}{\psi_j(\psigma)} 
		\tfrac{\psi_k(\sigma_k(\ozeta))}{\strut\lambda_k(\ozeta)} 
	}
	_{c^3_{kj}(\psigma,\ozeta)}
\d\ozeta
\end{align}
{and}
\begin{align}\label{coefftrafo}
\bm{a}(z, \zeta)
=  -\tfrac{1}{4}\lambda''_i(z) + \tfrac{1}{4} \lambda''_j(\zeta) + \frac{3(\lambda'_i(z))^2}{16\lambda_i(z)} 
-\frac{3(\lambda'_j(\zeta))^2}{16\lambda_j(\zeta)},
\end{align}%
if 
\begin{equation}\label{psiode}
 \psi'_i(\bm{\rho}) = \frac{\kappa_i({\phi_i^{-1}}(\bm{\rho}))}{2}\psi_i(\bm{\rho}), \quad \bm{\rho} \in (0,\phi_i(1)),
\end{equation}
$i=1,2,\ldots,n$,
holds. 
A solution of \eqref{psiode} 
is
\begin{equation}{\label{eq:psi}}
 \psi_i(\bm{\rho}) = \sqrt[4]{\frac{\lambda_i(\phi_i^{-1}(\bm{\rho}))}{\lambda_i(0)}}
\end{equation}
so that the transformation \eqref{secondcc} is invertible. Note that \eqref{eq:psi} also yields $\psi_j(\psigma)$ by applying the corresponding substitutions.

After utilizing the previous results, lengthy but straightforward algebraic manipulations yield the BCs 
\begin{subequations}\label{canBC}
\begin{align}
&  [\widetilde{\pmb{K}}(\bm{\rho},0) \tikzmark[3pt]{gleich}{$=$} 0]_{\lambda_i \geq \lambda_j,\, j \leq m}\\	 
&	[\widetilde{\pmb{K}}_{\bm{\sigma}}(\bm{\rho},0) 
 + (\underbrace{\,
		\textstyle \frac{\lambda'_j(0)}{4\sqrt{\lambda_j(0)}} 
		+ q_j\sqrt{\lambda_j(0)}\,
	}
	_{c^4_{j}})
	\widetilde{\pmb{K}}(\bm{\rho},0)\nonumber\\
	\tikzmark[0pt]{line}{} \strut
\tikz[overlay, remember picture]{
\node at (gleich.west|-line)[anchor=base west, inner sep=3pt]
{$= \widetilde{\pmb{\mathcal{C}}}[\widetilde{K}](\bm{\rho})]_{\lambda_i \geq \lambda_j,\, j > m}$};
}
\end{align}		
of the kernel PIDE \eqref{prenewtrans4} 
with
\begin{align}\label{eq:C_til}
&\widetilde{\pmb{\mathcal{C}}}[\widetilde{K}](\bm{\rho}) = 
	\underbrace{
		\tfrac{\sqrt{\lambda_j(0)}}{\psi_i(\prho)}\pA_0(z_i(\prho))
	}
	_{\bm{c}^5(\prho)} \\ \notag
& -\int_0^{z_i(\prho)} 
	\sum\limits_{k=1}^n 
	\widetilde{K}_{ik}(\prho,\sigma_k(\ozeta))
	\underbrace{
		A_{0,kj}(\ozeta)
		\frac{\sqrt{\lambda_j(0)} \psi_k(\sigma_k(\ozeta))}{\strut\lambda_k(\ozeta)} 	
	}
	_{c^6_{kj}(\ozeta)}
	\d\ozeta
\end{align}
\end{subequations}
and for $i = j$
\begin{equation}\label{canBC3}
	\widetilde{\pmb{K}}(\bm{\rho},\bm{\rho}) 
	= \underbrace{ 
		-\frac{1}{2}\sqrt{\lambda_i(0)}
		\int_0 ^{\bm{\rho}}\big(\bm{A}({z_i(\orho)}) 
			+ \mu_c\big)\d \orho
	}
	_{c^7_{i}(\bm{\rho})}.
	\end{equation}
Additionally, for $i \neq j$
\begin{subequations}\label{canBC2}
\begin{equation}
\widetilde{\pmb{K}}
\big(\bm{\rho},\sigma{_j}(z{_i}(\prho))\big) = 0
\end{equation}
and
\begin{equation}
\widetilde{\pmb{K}}_{\bm{\rho}}\big(\bm{\rho},\sigma{_j}(z{_i}(\prho))\big)  
=\underbrace{
	\textstyle
	\frac{
		\sqrt[4]{\rule[0.8ex]{0pt}{1ex}\lambda_i(0)\lambda_j(0)\lambda_i(z{_i}(\bm{\rho}))
			\lambda^3_j(z{_i}(\bm{\rho}))
		}
	}{
		\rule[0.8ex]{0pt}{1ex}\lambda_j(z{_i}(\bm{\rho})) - \lambda_i(z{_i}(\bm{\rho}))
	}
	\bm{A}(z{_i}(\bm{\rho}))
}
_{\bm{c}^8(\bm{\rho})}
\end{equation}
\end{subequations}
hold.

\subsubsection{Canonical kernel equations}
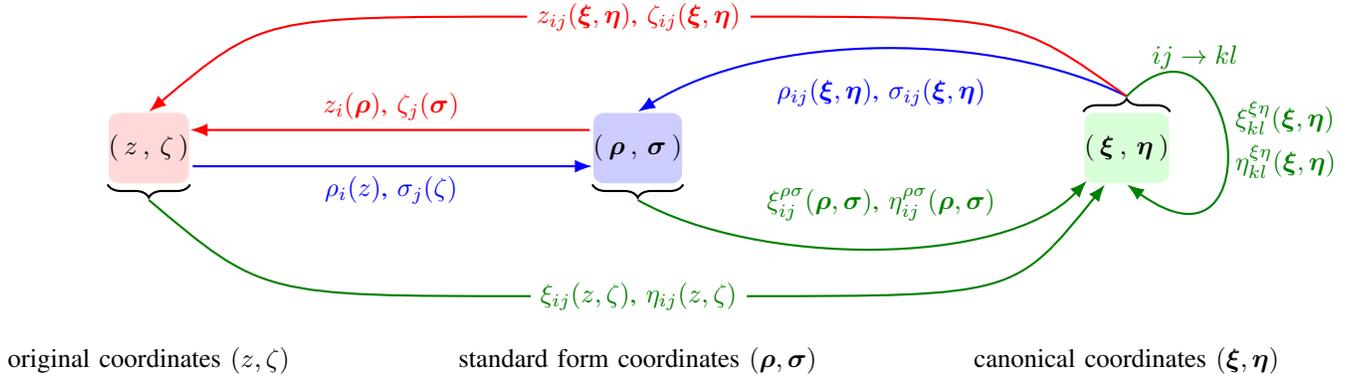
\begin{figure*}
\centering
\newcommand{\vertdist}{0cm}
\newcommand{\hordistbig}{6.5cm}
\newcommand{\hordistK}{1cm}
\newcommand{\anga}{75}
\newcommand{\xshiftlabels}{1cm}
\newcommand{\yshiftlabels}{0.25cm}
\newcommand{\curveshift}{0.5cm}
\newcommand{\yshiftlabelsbig}{0.75cm}
\newcommand{\textshift}{2cm}

\newcommand{\braceshift}{0mm}
\newcommand{\braceamp}{5pt}
\newcommand{\controlshift}{5cm}
\newcommand{\radius}{2cm}

\colorlet{redcol}{red!15}
\colorlet{bluecol}{blue!20}
\colorlet{greencol}{green!15}

\colorlet{redfront}{red}
\colorlet{bluefront}{blue}
\colorlet{greenfront}{green!50!black}

\newcommand{\prhocol}{\textcolor{bluefront}{\prho}}
\newcommand{\psigmacol}{\textcolor{bluefront}{\psigma}}

\begin{tikzpicture}[remember picture,thick,>={Latex},innode/.style={inner xsep=2pt,inner ysep=4pt, outer ysep=-2pt},
fillr/.style= {fill=#1, rounded corners}]
\node[inner sep=0pt](K) at (0,0){};
\node(z_ges)[below=\vertdist of K,inner xsep=0pt,fillr=redcol]{\strut $(
\tikz[remember picture,baseline=(z.base)] \node[innode] (z) {\strut$z$};
,
\tikz[remember picture,baseline=(zeta.base)] \node[innode] (zeta) {\strut$\zeta$};
)$};
\draw[decorate,decoration={brace,amplitude = \braceamp}]([yshift=-\braceshift]z_ges.south east)--node[below=\braceamp,inner sep=0pt](underbracez){}([yshift=-\braceshift]z_ges.south west);

\node(K_bar)[right=\hordistbig of K.center,anchor=center, inner sep=0pt]{};

\node(rho_ges)[below=\vertdist of K_bar, inner xsep=0pt,fillr=bluecol]{\strut $(
\tikz[remember picture,baseline=(rho.base)] \node[innode] (rho) {\strut$\prho$};
,
\tikz[remember picture,baseline=(sigma.base)] \node[innode] (sigma) {\strut $\psigma$};
)$};

\node(G)[right=\hordistbig of K_bar.center,anchor=center, inner sep=0pt]{};
\node(xi_ges)[below=\vertdist of G, inner xsep=0pt,fillr=greencol]{\strut $(
\tikz[remember picture,baseline=(xi.base)] \node[innode] (xi) {\strut $\pxi$};
,
\tikz[remember picture,baseline=(eta.base)] \node[innode] (eta) {\strut $\peta$};
)$};

\node(zrho_ges)[fit=(z_ges)(rho_ges), inner sep=0pt]{};
\node(rhoxi_ges)[fit=(rho_ges)(xi_ges), inner sep=0pt]{};
\begin{scope}[color=blue]
	\draw[->, shorten <=1pt](z_ges.east south east) --node[below](rhoi){$\rho_i(z)$, $\sigma_j(\zeta)$} (rho_ges.west south west);
\end{scope}
\begin{scope}[color=red]
	\draw[->, shorten <=1pt](rho_ges.west north west) --node[above](zi){$z_i(\prho)$, $\zeta_j(\psigma)$} (z_ges.east north east);
\end{scope}
\draw[decorate,decoration={brace,amplitude = \braceamp}]([yshift=-\braceshift]rho_ges.south east)--node[below=\braceamp,inner sep=0pt](underbracerho){}([yshift=-\braceshift]rho_ges.south west);
\begin{scope}[color=green!50!black]
	\node(xiijsigma) at (rhoxi_ges.south)[yshift=-\yshiftlabels]{$\xi_{ij}^{\rho\sigma}(\prho,\psigma)$,
	$\eta_{ij}^{\rho\sigma}(\prho,\psigma)$};
	\draw[->,shorten >=-2pt](underbracerho.center) .. controls ([yshift=-\curveshift]xiijsigma.south west) and ([yshift=-\curveshift]xiijsigma.south east) .. (xi_ges.south west);
\end{scope}
\draw[decorate,decoration={brace,amplitude = \braceamp}]([yshift=\braceshift]xi_ges.north west)--node[above=\braceamp,inner sep=0pt](overbracexi){}([yshift=\braceshift]xi_ges.north east);
\begin{scope}[color=blue]
\node(rhoijxi) at (rhoxi_ges.north)[yshift=\yshiftlabels]{$\rho_{ij}(\pxi,\peta)$,
 $\sigma_{ij}(\pxi,\peta)$};
\draw[->](overbracexi.center) .. controls ([yshift=\curveshift]rhoijxi.north east) and ([yshift=\curveshift]rhoijxi.north west) .. (rho_ges.north);
\end{scope}
\node(zxi_ges)[fit=(z_ges)(xi_ges)(zi)(xiijsigma)]{};
\begin{scope}[color=green!50!black]
	\node(xiijz)[yshift=-\yshiftlabelsbig]at(zxi_ges.south){$\xi_{ij}(z,\zeta)$, $\eta_{ij}(z,\zeta)$};
	\draw[->,shorten >=0pt](underbracez.center) .. controls ([xshift=-\controlshift]xiijz) .. (xiijz) .. controls ([xshift=\controlshift]xiijz) .. (xi_ges.south south west);
\end{scope}

\begin{scope}[color=red]
\node(zijxi)[yshift=\yshiftlabelsbig]at(zxi_ges.north){$z_{ij}(\pxi,\peta)$, $\zeta_{ij}(\pxi,\peta)$};
\draw[->](overbracexi.center) .. controls ([xshift=\controlshift]zijxi) .. (zijxi) .. controls ([xshift=-\controlshift]zijxi) .. (z_ges.north);
\end{scope}

\node(phant1)[above right=\radius of overbracexi.north east,inner sep=0pt]{};
\node(phant2)[below right=\radius of xi_ges.south east,inner sep=0pt]{};
\begin{scope}[color=green!50!black]
\draw[->](overbracexi.center) .. node[above,near start]{$ij \rightarrow kl$}node[right, text width=1.5cm](xikl){$\xi_{kl}^{\xi\eta}(\pxi,\peta)$ \\[4pt] $\eta_{kl}^{\xi\eta}(\pxi,\peta)$}  controls (phant1) and (phant2) .. (xi_ges.south);
\end{scope}

\node(Ktext)[below=\textshift of z_ges]{\strut original coordinates $(z,\zeta)$};
\node(rhotext)[below=\textshift of rho_ges]{\strut standard form coordinates $(\prho,\psigma)$};
\node(xitext)[below=\textshift of xi_ges]{\strut canonical coordinates $(\pxi,\peta)$};

\node(mid1)[inner sep=0pt,outer sep=0pt]at($(K.east)!0.5!(K_bar.west)$){};
\node(mid2)[inner sep=0pt]at($(K_bar.east)!0.5!(G.west)$){};


\end{tikzpicture}
\caption{Overview of the changes of coordinates for mapping the kernel equations into their canonical form.}
\label{fig:coordinate_transformations}
\end{figure*}
The kernel PIDE \eqref{prenewtrans4} can now be mapped into the \emph{canonical form}. 
This only requires a linear change of coordinates so that no additional first-order partial derivatives are introduced in the result.
The calculation of the related transformation coincides with the method for classifying second-order PDEs (see, e.\:g., \cite[Ch. 4]{My07}). 

The change of coordinates mapping the kernel PIDE \eqref{prenewtrans4} for $\widetilde{\pmb{K}}(\bm{\rho},\bm{\sigma}) = \widetilde{K}_{ij}(\bm{\rho},\bm{\sigma})$ into its canonical form depends on the relation of the corresponding diffusion coefficients $\lambda_i$ and $\lambda_j$. In particular, one obtains
\begin{subequations}\label{nfcoord}
\begin{align}
 \lambda_i \geq \lambda_j:\ \  \bm{\xi} &= \xi{_{ij}^{\rho\sigma}}(\prho,\psigma) = \bm{\rho} + \bm{\sigma}\\
  \bm{\eta} &= \eta{_{ij}^{\rho\sigma}}(\prho,\psigma) = \bm{\rho} - \bm{\sigma}\label{cononcc1}\\
\lambda_i < \lambda_j:\ \ \bm{\xi}  &= \xi{_{ij}^{\rho\sigma}}(\prho,\psigma) = \phi_i(1) + \phi_j(1) - (\bm{\rho} + \bm{\sigma})\label{xil}\\
\bm{\eta} &= \eta{_{ij}^{\rho\sigma}}(\prho,\psigma) = -(\phi_i(1) - \phi_j(1)) + \bm{\rho} - \bm{\sigma}\label{etal}
\end{align}
\end{subequations}
with $\phi_i$ defined in \eqref{phidef}.
In \eqref{nfcoord}, $(\pxi,\peta)$ denotes a point in the canonical coordinate system, whereas $\xi_{ij}^{\rho\sigma}(\prho,\psigma)$ and $\eta_{ij}^{\rho\sigma}(\prho, \psigma)$ are the corresponding transformations. Thereby, the superscript ${\rho\sigma}$ is introduced to identify the considered original coordinates.
Besides mapping the kernel PIDE in its canonical form, the coordinates \eqref{xil} and \eqref{etal} additionally lead to spatial domains for $\lambda_i < \lambda_j$, that coincide with the spatial domains in the case $\lambda_i \geq \lambda_j$. This simplifies the subsequent derivation of the kernel integral equations. By adding and subtracting the equations for $\bm{\xi}$ and $\bm{\eta}$ in \eqref{nfcoord}, the inverse change of coordinates
\begin{subequations}\label{nfcoordi}
	\begin{align}
	\lambda_i \geq \lambda_j:\quad \bm{\rho}  &= \rho{_{ij}}(\pxi,\peta) =  \frac{\bm{\xi} + \bm{\eta}}{2} \\
	\bm{\sigma} &= \sigma{_{ij}}(\pxi,\peta) =  \frac{\bm{\xi} -\bm{\eta}}{2} \\
	\lambda_i < \lambda_j:\quad \bm{\rho}  &= \rho{_{ij}}(\pxi,\peta) = \phi_i(1) - \frac{\bm{\xi} - \bm{\eta}}{2} \\ 
	\quad \bm{\sigma} &=  \sigma{_{ij}}(\pxi,\peta) = \phi_j(1) - \frac{\bm{\xi} + \bm{\eta}}{2}\label{etali}
	\end{align}
\end{subequations}
is readily found. 
Note that a point $(\prho,\psigma)$ can either be computed by the functions $\rho_i(z)$, $\sigma_j(\zeta)$ (see \eqref{cchange}) from the original coordinates, or by the functions $\rho_{ij}(\pxi,\peta)$ and $\sigma_{ij}(\pxi,\peta)$ according to \eqref{nfcoordi} from the canonical coordinates.
In order to simplify the notation, introduce
\begin{equation}\label{shdef}
	\bm{s} = s_{ij} = \begin{cases}
	1, & \lambda_i \geq \lambda_j\\
	-1, & \lambda_i < \lambda_j.
	\end{cases}
\end{equation}
Insert \eqref{cchange} in \eqref{nfcoord} and use \eqref{shdef} to obtain the \emph{overall nonlinear coordinate transformation}
\begin{subequations}\label{eq:one_step_trafo}
\begin{align}
	\pxi &= \xi_{ij}(z,\zeta) = \tfrac{1}{2}(1-\ps)(\phi_i(1) + \phi_j(1)) + \ps(\phi_i(z) + \phi_j(\zeta)) \\
	\peta &=\eta_{ij}(z,\zeta) = -\tfrac{1}{2}(1-\ps)(\phi_i(1) - \phi_j(1)) + \phi_i(z) - \phi_j(\zeta)
\end{align}
\end{subequations}
from the original $(z,\zeta)$-coordinates into the canonical coordinates.
This transformation can be solved for $(z,\zeta)$, which yields
\begin{subequations}\label{eq:inverse_one_step_trafo}
	\begin{align}
			z &= z_{ij}(\pxi,\peta) =  \phi_i^{-1}(\tfrac{1}{2}(\ps\pxi+\peta)+\tfrac{1}{2}(1-\ps)\phi_i(1)) \\
			\zeta &= \zeta_{ij}(\pxi,\peta) = \phi_j^{-1}(\tfrac{1}{2}(\ps\pxi-\peta)+\tfrac{1}{2}(1-\ps)\phi_j(1)).
			\label{eq:zeta_inv}
	\end{align}
\end{subequations}
By combining the transformations \eqref{eq:one_step_trafo} and \eqref{eq:inverse_one_step_trafo}, 
the canonical coordinates of different indices $i$ and $j$ can always be transformed into each other, which is necessary due to the coupling between different kernel PIDEs. In particular, the transformations%
\begin{subequations}\label{eq:trafo_xi_eta}%
	\begin{align}%
		&\xi_{kl}^{\xi\eta}(\pxi,\peta) =
		\xi_{kl}\Big(z_{ij}(\pxi,\peta),\zeta_{ij}(\pxi,\peta)\Big) \\
		&\eta_{kl}^{\xi\eta}(\pxi,\peta) =
		\eta_{kl}\Big(z_{ij}(\pxi,\peta),\zeta_{ij}(\pxi,\peta)\Big)
	\end{align}
\end{subequations}
can be introduced to transform from one canonical coordinate system $(\pxi,\peta) = (\xi_{ij},\eta_{ij})$ into the canonical coordinate system of a different element $(\xi_{kl}, \eta_{kl})$.
Observe that the canonical coordinates can now be calculated from the original coordinates by \eqref{eq:one_step_trafo}, from the $(\prho,\psigma)$-coordinates by \eqref{nfcoord} or from the canonical coordinates of a different kernel element by \eqref{eq:trafo_xi_eta}.
An overview of the various introduced coordinate transformations is given in Figure \ref{fig:coordinate_transformations}.

Now define $\bm{G}(\bm{\xi},\bm{\eta}) = \widetilde{\pmb{K}}(\rho{_{ij}}(\pxi,\peta),\sigma{_{ij}}(\pxi,\peta))$, {which can be expressed as $\bm{G}(\xi_{ij}^{\rho\sigma}(\prho,\psigma),\eta_{ij}^{\rho\sigma}(\prho,\psigma)) = \widetilde{\pmb{K}}(\prho,\psigma)$ in the standard form coordinates, differentiate it twice \wrt $\prho$ and $\psigma$ and insert it with} 
\eqref{nfcoord} {in} \eqref{prenewtrans4}. {Using \eqref{shdef} provides the compact notation of} the \emph{canonical kernel PIDE}
\begin{equation}\label{canonicalpde}
 4\bm{s}\,\bm{G}_{\bm{\xi}\bm{\eta}}(\bm{\xi},\bm{\eta}) =   \breve{\bm{\mathcal{B}}}[G](\bm{\xi},\bm{\eta})
\end{equation}
after straightforward intermediate computations. %
Thereby, the abbreviation
\begin{align}\notag
& \breve{\bm{\mathcal{B}}}[G](\pxi,\peta) =
\Big(
	\pa\big(z_{ij}(\pxi,\peta),\zeta_{ij}(\pxi,\peta) \big) + \mu_c 
\Big)\pG (\pxi,\peta) \\  \notag
& - \bm{c}^{1}(\rho_{ij}(\pxi,\peta),\sigma_{ij}(\pxi,\peta)) \\\notag
&+ \sum\limits_{k=1}^{n} \bigg( 
	G_{ik}(\xi^{\xi\eta}_{ik}(\pxi,\peta),\eta^{\xi\eta}_{ik}(\pxi,\peta))  c^2_{kj}(\sigma_{ij}(\pxi,\peta))  \\ \notag
	&  + \int_{\zeta_{ij}(\pxi,\peta)}^{z_{ij}(\pxi,\peta)}
		G_{ik}
		\Big(
			\xi_{ik}\big(
				z_{ij}(\pxi,\peta),
				\ozeta
			\big),
			\eta_{ik}\big(
				z_{ij}(\pxi,\peta),
				\ozeta
			\big)
		\Big) \\ \label{eq:B_breve}
		& \hspace{1.8cm} \cdot c^3_{kj}(\sigma_{ij}(\pxi,\peta),\ozeta) 
	\d\ozeta
\bigg)
\end{align}
was used (see \eqref{eq:B}).
\begin{rem}
 It should be noted that the proposed three step transformation into the kernel PIDE \eqref{canonicalpde} is much simpler as directly applying the classical change of coordinates for classifying second-order PDEs to the kernel PIDEs \eqref{keqcompdiagpde} and \eqref{keqcomp1pde}. 
This is due to the fact that the latter approach yields first-order partial derivatives with coefficient functions which exhibit an involved structure.
Consequently, the determination of the corresponding transformation for their elimination is impeded.
 \hfill $\triangleleft$
\end{rem}
\begin{rem}
 For mutually different \emph{constant diffusion coefficients} $\lambda_{i} = const.$, $i = 1,2,\ldots,n$, the overall linear change of coordinates 
 \begin{subequations}
 \begin{align}
   \lambda_i \geq \lambda_j: \quad \bm{\xi} &= \textstyle\frac{z}{\sqrt{\lambda_i}} + \textstyle\frac{\zeta}{\sqrt{\lambda_j}} \quad \text{and} \quad \bm{\eta} = \textstyle\frac{z}{\sqrt{\lambda_i}} - \textstyle\frac{\zeta}{\sqrt{\lambda_j}}\\
   \lambda_i < \lambda_j:\quad \bm{\xi}  &= \textstyle\frac{1}{\sqrt{\lambda_i}} + \textstyle\frac{1}{\sqrt{\lambda_j}} - (\frac{z}{\sqrt{\lambda_i}} + \textstyle\frac{\zeta}{\sqrt{\lambda_j}})\\
	\bm{\eta} &= -(\textstyle\frac{1}{\sqrt{\lambda_i}} + \textstyle\frac{1}{\sqrt{\lambda_j}}) + \textstyle\frac{z}{\sqrt{\lambda_i}} - \textstyle\frac{\zeta}{\sqrt{\lambda_j}}
	\end{align}
 \end{subequations}
 can be applied to obtain the canonical kernel PIDE in a single step. This directly follows from \eqref{cchange}, \eqref{phidef} and \eqref{nfcoord}. \hfill $\triangleleft$
\end{rem}

For the derivation of the canonical kernel equations it remains to determine the BCs for \eqref{canonicalpde}. 
They
result from applying \eqref{nfcoord} and \eqref{nfcoordi} to \eqref{canBC}--\eqref{canBC2}. For this, the spatial domain of \eqref{canonicalpde} is investigated.

\begin{figure}[t]
	\centering	
\usetikzlibrary{intersections}

\newcommand{\coordsys}[2][1]{
\begin{tikzpicture}[scale=#1]
\newcommand{\linew}{2pt} 
\newcommand{\alphaij}{#2}
\newcommand{\beschbreit}{0.025}
\newcommand{\achsenfaktor}{1.4}
\newcommand{\achsenfaktorx}{1.2}
\draw[->](0,0)--(\achsenfaktorx+\achsenfaktorx*\alphaij,0)node[below right]{$\pxi$};
\draw[->](0,-\achsenfaktor+\achsenfaktor*\alphaij)--(0,\achsenfaktor*\alphaij)node[left]{$\peta$};

\coordinate(a) at (\alphaij,0);
\draw[shift={(a)}](0,-\beschbreit)--(0,+\beschbreit)node[above]{$\;\;\,a_{\phantom{ij}}$};	
\coordinate(1pa) at (1+\alphaij,0);
\draw[shift={(1pa)}](0,-\beschbreit)--(0,+\beschbreit)node[above]{$\phi_i(1)+\phi_j(1)$};	

\coordinate(b) at (0,-1+\alphaij);
\draw[shift={(b)}](\beschbreit,0)--(-\beschbreit,0)node[left]{$b$};	
\coordinate(a2) at (0,\alphaij);
\draw[shift={(a2)}](\beschbreit,0)--(-\beschbreit,0)node[left]{$a$};
\coordinate(0) at (0,0);
\draw[shift={(0)}](\beschbreit,0)--(-\beschbreit,0)node[left]{$0$};

\node at (\alphaij+0.55,\alphaij/2.5-0.4){\normalsize$\mathcal{D}_{ij}$};

\draw[line width=\linew](0,0)--node[pos=0.25,xshift=-0.1cm](gam1){}node[sloped,pos=0.4,above]{$\Gamma_{\!1}$}(\alphaij,\alphaij);
\draw[dashed](\alphaij,\alphaij)--(1+\alphaij,-1+\alphaij);
\coordinate(newx) at (1+\alphaij,-1+\alphaij);

\begin{scope}[x={(newx)},y={($(0,0)!1!90:(newx)$)}]
\draw[out=-40,in=170, line width=\linew,name path=gam2path](0,0) to (0.5,0);
\draw[out=170-180,in=175, line width=\linew](0.5,0) to (1,0);
\draw[out=-40,in=170, line width=\linew, shift={(0,-0.03)}](0-0.01,0) to node[sloped,pos=0.58,below]{$\Gamma_{\!2}: \peta_{l}(\pxi)$}node[inner sep=0pt,pos=1, outer sep=0pt](gam2){}(0.5,0);
\draw[out=170-180,in=175, line width=\linew, shift={(0,-0.03)}](0.5,0) to (1+0.005,0);
\end{scope}
\path[name path=gam1proj](gam1)--(gam1|-b);

\node (gam2del) at ([yshift=0.42cm]gam2){};
\draw[line width=3.5pt,white,shorten <=1pt](gam2)--(gam2del); 
\draw[line width=0.5pt,->](gam2)--(gam2del);

\begin{scope}[scale=1/#1] 
\path[name intersections={of=gam2path and gam1proj, by=schnitt}];
\end{scope}
\node (gam2aaadel) at ([xshift=0.33cm]schnitt){};
\draw[densely dash dot,line width=0.5pt,->](schnitt)--(gam2aaadel);

\node (gam1del) at ([xshift=0.33cm]gam1){};
\draw[densely dash dot,line width=0.5pt,->](gam1)--(gam1del);

\end{tikzpicture}

}
	{\small 
		\coordsys[2.75]{0.6} 
	}
	\caption{Spatial domains $\mathcal{D}_{ij}$ of the canonical kernel equations \eqref{canonkeq} with
		$\lambda_i \geq \lambda_j:\; a = \phi_i(1), b = \phi_i(1) -\phi_j(1)$ and
		$\lambda_i < \lambda_j:\; a = \phi_j(1), b = -(\phi_i(1) -\phi_j(1))$.
		The bold lines at $\Gamma_{\!1}$ and $\Gamma_{\!2}$ represent the BCs of the kernel BVP, in which $\Gamma_{\!1}$ is an artificial BC for $\lambda_i < \lambda_j$ (see \eqref{wpbc}). The dash dotted arrows ('$-\cdot-$') indicate the integration in the $\pxi$-direction and the solid arrow ('--') the integration in the $\peta$-direction. Note that in the case $i=j$ the 
		spatial domain is the triangle $\peta \in (0,\phi_i(1))$, $\pxi \in (\peta,2\phi_i(1)-\peta)$, as $b=0$ and $\peta_l(\pxi) \equiv 0$.}\label{Fig1}
\end{figure}
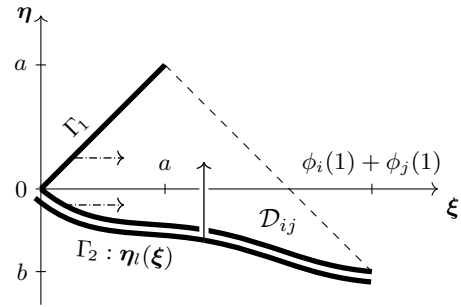

Consider the boundary $\zeta = z$, $z \in [0,1]$, of the spatial domain \wrt the original kernel PIDEs \eqref{keqcompdiagpde} and \eqref{keqcomp1pde} and substitute $\zeta=z$ in \eqref{eq:one_step_trafo} to obtain
\begin{subequations}
	\begin{align}
		\pxi &= \tfrac{1}{2}(1-\ps)(\phi_i(1) + \phi_j(1)) + \ps\bm{\beta}(z) \label{xif}\\
		\peta &= -\tfrac{1}{2}(1-\ps)(\phi_i(1) - \phi_j(1)) + \phi_i(z) - \phi_j(z) \label{etaf}
	\end{align}
\end{subequations}
with $\bm{\beta}(z) = \phi_i(z) + \phi_j(z)$.
The function $\bm{\beta}(z)$
is strictly monotonically increasing, because $\phi_i$ and $\phi_j$ have this property (see \eqref{phidef}). Hence, the inverse $\bm{\beta}^{-1}$ exists and thus \eqref{xif} can be solved for $z$, resulting in
\begin{align}
	\pz_{l}(\pxi) = {\bm{\beta}}^{-1}\Big(\ps\pxi+\tfrac{1}{2}(1-\ps)(\phi_i(1)+\phi_j(1))\Big).
\end{align}
Inserting this into \eqref{etaf} yields the lower boundary
\begin{equation}\label{eq:eta_l}%
\peta_{l}(\pxi) = -\tfrac{1}{2} (1-\ps)(\phi_i(1)-\phi_j(1)) + \phi_i(\pz_{l}(\pxi)) - \phi_j(\pz_{l}(\pxi))%
\end{equation}%
of the domain $\mathcal{D}_{ij}$ in Figure \ref{Fig1}.
It is easy to show that {for $i\neq j$,} $\peta_l'(\pxi) < 0$ for $\bm{\xi} \in [0,\phi_i(1)+\phi_j(1)]$ so that $\bm{\eta}_l(\bm{\xi})$ is strictly monotonically decreasing. 
This implies that the resulting boundary of \eqref{canonicalpde} only evolves in the fourth quadrant of the $(\bm{\xi},\bm{\eta})$-plane (see Figure \ref{Fig1}). In contrast,{ inserting $i=j$ in \eqref{eq:eta_l} shows that}
\begin{align}\label{eq:eta_l_0}%
	[\peta_l(\pxi) = 0]_{i=j}
\end{align}
{holds} so that the spatial domain is the triangle $\peta \in (0,\phi_i(1))$, $\pxi \in (\peta,2\phi_i(1)-\peta)$.
By a similar but very elementary reasoning one can transform 
the remaining boundaries of the kernel PIDEs \eqref{keqcompdiagpde} and \eqref{keqcomp1pde} into the $(\bm{\xi},\bm{\eta})$-coordinates with \eqref{eq:one_step_trafo}. 
The resulting spatial domain $\mathcal{D}_{ij}$ is depicted in Figure \ref{Fig1}.

With the spatial domain of the canonical coordinates known, the BCs \eqref{canBC}--\eqref{canBC2} can be mapped into $(\bm{\xi},\bm{\eta})$-coordinates by making use of \eqref{nfcoord} and \eqref{nfcoordi}, resulting in the \emph{canonical kernel equations}
\begin{subequations}\label{canonkeq}
	\newcommand{\spleft}{-1cm}
	\begin{flalign}
	  &\hspace{\spleft}&4\bm{s}\,\bm{G}_{\bm{\xi}\bm{\eta}}(\bm{\xi},\bm{\eta}) &=
	  \breve{\bm{\mathcal{B}}}[G](\bm{\xi},\bm{\eta})\label{ckeqpde} \hspace{3.2cm}\\
	   &\hspace{\spleft}&[\bm{G}(\bm{\eta},\bm{\eta}) &= 0]_{\lambda_i \geq \lambda_j,\, j \leq m} \label{eq:canonBC1}\\
	   &\mathrlap{[\bm{G}_{\bm{\xi}}(\bm{\eta},\bm{\eta}) - \bm{G}_{\bm{\eta}}(\bm{\eta},\bm{\eta})   + c^4_{j}\bm{G}(\bm{\eta},\bm{\eta}) = \breve{\bm{\mathcal{C}}}[G](\bm{\eta})]\rightcond{\lambda_i&\geq\!\lambda_j\\ j&>\!m}} \label{eq:canonBC2}\\
	   &\hspace{\spleft}&[\bm{G}(\bm{\xi},0) &= c^7_{i}(\tfrac{\bm{\xi}}{2})]_{i = j} \label{eq:canonBC4}\\
	   &\hspace{\spleft}&[\pG(\peta,\peta) &= \bm{g}_f(\peta)]_{\lambda_i<\lambda_j} \label{wpbc} \\	   
	  &\hspace{\spleft}&[\bm{G}(\bm{\xi},\bm{\eta}_l(\bm{\xi})) &= 0]_{i \neq j} \label{eq:canonBC5}\\
	  &\mathrlap{[\ps\bm{G}_{\bm{\xi}}(\bm{\xi},\bm{\eta}_l(\bm{\xi})) +  \bm{G}_{\bm{\eta}}(\bm{\xi},\bm{\eta}_l(\bm{\xi})) = \bm{c}^8\big(\rho_{ij}(\pxi,\peta_l(\pxi))\big)]_{i \neq j}} \label{eq:canonBC6}
	\end{flalign} 
\end{subequations}
for $i,j = 1,2,\ldots,n$, in which \eqref{ckeqpde} is defined on $\mathcal{D}_{ij}$ (see Figure \ref{Fig1}) and
\begin{align}\label{eq:C_breve}
&\breve{\bm{\mathcal{C}}}[G](\bm{\eta}) =
\bm{c}^5(\rho_{ij}(\peta,\peta))  \\\notag
	& - \int\limits_0^{\mathclap{z_{ij}(\peta,\peta)}} 
	\hspace{0.8em}  \sum\limits_{k=1}^n 
 G_{ik}
	\Big(
	\xi_{ik}\big(z_{ij}(\peta,\peta),\ozeta\big),\eta_{ik}\big(z_{ij}(\peta,\peta),\ozeta\big)
	\Big) 	c^6_{kj}(\ozeta)  \d\ozeta
\end{align}
was used (see \eqref{eq:C_til}).
%
%
The BC \eqref{wpbc} at $\Gamma_{\!1}$ is an \emph{artificial BC}, that has to be introduced in order to ensure well-posed kernel equations.  Thereby, $\bm{g}_f \in C^2[0,\phi_j(1)]$ is a degree of freedom and can be chosen arbitrarily.
\begin{rem}
 It is noteworthy that two BCs are needed at the boundary $\Gamma_{\!2}$, because the derivation of the integral equations requires one BC at $\Gamma_{\!2}$ for the integration in the direction of $\pxi$ and {one BC at $\Gamma_{\!2}$ for the integration in the direction of} $\peta$ (see Figure \ref{Fig1}). 
 Furthermore, the  BC{s} at $\Gamma_{\!1}$ and  $\Gamma_{\!2}$ utilized for the integration \wrt the $\pxi$-direction imply that the $\pxi$-axis within $\mathcal{D}_{ij}$ constitutes a \emph{line of separation} defining two different parts of the corresponding kernel on {the} spatial domain $\mathcal{D}_{ij}${,} leading to piecewise defined kernels for $i \neq j$. It can be shown that the line of separation is no longer a straight line in the original coordinates but a strictly monotonically increasing curve defined by $\zeta = \phi_j^{-1}(\phi_i(z))$, $z \in [0,1]$, for $\lambda_i > \lambda_j$ and $z = \phi_i^{-1}(\phi_i(1) - \phi_j(1) + \phi_j(\zeta))$, $\zeta \in [0,1]$, if $\lambda_i < \lambda_j$. \hfill $\triangleleft$
\end{rem}

\subsection{Kernel Integral Equations}\label{sec:kieq}
In order to solve the canonical kernel equations \eqref{canonkeq} utilizing a successive approximation, they are transformed into integral equations. The convergence analysis can be simplified by rewriting the second-order PIDE \eqref{ckeqpde} into a system of two first-order PIDEs. By introducing
\begin{equation}\label{eq:Hdef}
\pH(\pxi,\peta) = \pG_{\pxi}(\pxi,\peta),
\end{equation}
the PIDE \eqref{ckeqpde} can be written as 
\begin{subequations}\label{eq:firstorderkeq}
\begin{align} 
 \pG_{\pxi}(\pxi,\peta) &= \pH(\pxi,\peta) \label{eq:G_xi} \\
\pH_{\peta}(\pxi,\peta) &= \tfrac{1}{4\ps}\breve{\pmb{\mathcal{B}}}[G](\pmb{\xi},\pmb{\eta}). \label{eq:H_eta}
\end{align}
\end{subequations}
To formulate the corresponding BCs, 
solve \eqref{eq:canonBC2} for $\pG_{\peta}(\peta,\peta)$, insert the result into
\begin{equation}
 \d_{\peta} \pG(\peta,\peta) = \pG_{\pxi}(\peta,\peta) + \pG_{\peta}(\peta,\peta)
\end{equation}
and use \eqref{eq:Hdef}.
Then, an integration \wrt $\peta$ yields
\begin{flalign} \label{eq:G_eta_eta}
	\Big[
	\pG(\peta,\peta) =\!\! \int\limits_{0}^{\peta}\!\! \Big( 2\pH(\oeta,\oeta) + c^4_{j} \pG(\oeta,\oeta) - \breve{\pmb{\mathcal{C}}}[G](\oeta)\Big)\d\oeta
	\Big]\rightcond{\lambda_i&\geq\!\lambda_j\\ j&>\!m}.
\end{flalign}
In view of \eqref{eq:canonBC5} and \eqref{eq:Hdef} one obtains
\begin{align}\label{eq:canonBC5_diff}
	\d_{\pxi} \pG(\pxi,\peta_l(\pxi)) 
	= \pH(\pxi,\peta_l(\pxi)) + \pG_{\peta}(\pxi,\peta_l(\pxi)) \, \peta^{\prime}_l(\pxi) = 0.
\end{align}
Solving \eqref{eq:canonBC6} for $\pG_{\peta}(\pxi,\peta_l(\pxi))$ and inserting the result in \eqref{eq:canonBC5_diff} yields
\begin{equation}\label{eq:Hxietal}
\Big[\pH(\pxi,\peta_l(\pxi)) = \underbrace{\frac{\pmb{c}^8\big(\rho_{ij}(\pxi,\peta_l(\pxi))\big) \: \peta_l'(\pxi)}{\ps\peta_l'(\pxi)-1}}_{\pmb{c}^9(\pxi)}\Big]_{i\neq j}
\end{equation}
after a simple rearrangement.
It can easily be shown that $-1 < \peta_l'(\pxi) < 0$ for $\pxi \in [0, \phi_i(1) + \phi_j(1)]$. Hence, the denominator in \eqref{eq:Hxietal} cannot be zero. Finally, differentiate \eqref{eq:canonBC4} \wrt $\pxi$ and use \eqref{eq:Hdef}, the definition of $c_i^7$ in \eqref{canBC3} and \eqref{eq:eta_l_0} to obtain
\begin{equation} \label{eq:H_xi_0}
 \big[
 \pH(\pxi,0) = \pH(\pxi,\peta_l(\pxi)) = \underbrace{-\tfrac{\sqrt{\lambda_i(0)}}{4}\big(\pmb{A}(z_{ij}(\pxi,0)) +\mu_c
 \big)}_{\pmb{c}^{10}(\pxi)}
 \big]_{i=j}.
\end{equation}
For the ease of presentation, the left boundary of the spatial domain $\mathcal{D}_{ij}$ can be introduced as 
\begin{equation}\label{eq:xi_l}
 \pxi_l(\peta) = \begin{cases}
  \peta, & \peta \geq 0 \\
  \peta_l^{-1}(\peta), & \peta < 0
 \end{cases}
\end{equation}
with $\peta_l^{-1}$ being the inverse of \eqref{eq:eta_l}.
With this, an alternative way to write \eqref{eq:canonBC5} is 
\begin{equation}\label{eq:canonBC5_alt}
[\pmb{G}(\pxi_l(\peta),\peta) = 0]_{i \neq j,\, \peta<0}.
\end{equation}
Now, the BC on the left boundary of the spatial domain $\mathcal{D}_{ij}$ can be compactly written as
\begin{align}
	\pG(\pxi_l(\peta),\peta) = 
	\begin{cases}
		\pG(\peta,\peta), & \peta \geq 0\\
		0, & i \neq j,\ \peta < 0
	\end{cases}
\end{align}
in view of \eqref{eq:xi_l} and \eqref{eq:canonBC5_alt}.
With this and $\pG(\peta,\peta)$ determined by \eqref{eq:canonBC1}, \eqref{wpbc} and \eqref{eq:G_eta_eta}, the BCs needed for the derivation of the integral equations are
\begin{subequations}
	\begin{align}
		\pG(\pxi_l(\peta),\peta) &= \label{eq:G_BC}
		\begin{cases}	
			\int_{0}^{\peta}\big( 2\pH(\oeta,\oeta) + \mathrlap{c^4_{j} \pG(\oeta,\oeta)  
			 - \breve{\pmb{\mathcal{C}}}[G](\oeta) \big)\d\oeta,} \\
			  & \lambda_i \geq \lambda_j,\ j>m,\ \peta\geq0 \\[4pt]
		0, &  \lambda_i \geq \lambda_j,\ j \leq m,\ \peta\geq 0 \\
		\pmb{g}_f(\peta) ,&\lambda_i<\lambda_j,\ \hphantom{ j \leq m,}\ \, \tikzmark[0pt]{peta>0}{$\peta\geq 0$} \\
		0, & \tikzmark[0pt]{ineqj}{$i \neq j$},
		\end{cases}
		\tikz[overlay,remember picture]{
		\node at (ineqj-|peta>0){$\peta<0$};
		} \\
\intertext{and}
		\pH(\pxi,\peta_l(\pxi)) &= \begin{cases} \label{eq:H_BC}
			\pmb{c}^{10}(\pxi), & i=j \\
			\pmb{c}^{9}(\pxi), & i\neq j
		\end{cases}.
	\end{align}
\end{subequations}

Now, \eqref{eq:G_xi} is integrated \wrt $\pxi$ and $\eqref{eq:H_eta}$ \wrt $\peta$, leading to 
\begin{subequations}
	\begin{align}
		 \pG(\pxi,\peta) &=\pG(\pxi_l(\peta),\peta) + \int_{\pxi_l(\peta)}^{\pxi} \pH(\oxi,\peta) \d \oxi  \label{eq:inteq1}\\
		 \pH(\pxi,\peta) &= \pH(\pxi,\peta_l(\pxi)) + \int_{\peta_l(\pxi)}^{\peta}   \tfrac{1}{4\ps}\breve{\bm{\mathcal{B}}}[G](\pmb{\xi},\oeta) \d\oeta, \label{eq:inteq2}
	\end{align}
\end{subequations}
when taking the boundaries of $\mathcal{D}_{ij}$ into account.
After inserting \eqref{eq:G_BC} with \eqref{eq:C_breve} into \eqref{eq:inteq1}, 
the first \emph{kernel integral equation} reads
\begin{align}\label{eq:inteq1_norm}%
	\pG(\pxi,\peta) &= \pG_0(\peta) + \pmb{F}_{G}[G,\pH](\pxi,\peta) = \widetilde{\pmb{F}}_{G}[G,\pH](\pxi,\peta),
\end{align}%
where
\begin{align}\label{eq:G0}%
		\pG_0(\peta) &= 
		\begin{cases}
			-\int_0^{\peta} \pmb{c}^5(\rho_{ij}(\oeta,\oeta)) \d\bar{\eta}, & \lambda_i \geq \lambda_j,\ j>m,\ \peta \geq 0  \\
			\pmb{g}_f(\peta), &\lambda_i < \lambda_j,\ \hphantom{j>m,\ \,} \peta \geq 0  \\
			0, & \text{else}
		\end{cases}
	\end{align}%
and
\begin{align}\label{eq:int_operatorG}%
	&\pmb{F}_G[G,\pH](\pxi,\peta) =\int_{\pxi_l(\peta)}^{\pxi} \pH(\oxi,\peta)\d\oxi \\ \notag
	& + \Big[
	\int_0^{\peta}
	\Big(
		2 \pH (\oeta,\oeta) + c^4_{j} \pG(\oeta,\oeta) 
		+ 
		\int_0^{z_{ij}(\oeta,\oeta)}
			\sum_{k=1}^{n}
				c^6_{kj}(\ozeta)\\ \notag
				& \cdot G_{ik}
				\big(
					\xi_{ik}(
						z_{ij}(\oeta,\oeta),\ozeta
					)
					,\eta_{ik}(
						z_{ij}(\oeta,\oeta),\ozeta
					)
				\big) \d\ozeta 
	\Big)\d\oeta
	\Big]\rightcond{\lambda_i&\geq\!\lambda_j\\ j&>\!m\\\peta &>\!0}
\end{align}%
hold.

Furthermore, after inserting \eqref{eq:B_breve} in \eqref{eq:inteq2} 
the second \emph{kernel integral equation} follows as
\begin{align}\label{eq:inteq2_norm}
	\pH(\pxi,\peta) &= \pH_0(\pxi) + \pmb{F}_{H}[G](\pxi,\peta) = \widetilde{\pmb{F}}_{H}[G](\pxi,\peta)%
\end{align}
with 
\begin{align}\label{eq:H0}
	\pH_0(\pxi) &=
	\pH(\pxi,\peta_l(\pxi))
	- \frac{1}{4\ps}\int_{\peta_l(\pxi)}^{\peta} \pmb{c}^1(\rho_{ij}(\pxi,\oeta),\sigma_{ij}(\pxi,\oeta))\d\oeta
\end{align}%
and $\pH(\pxi,\peta_l(\pxi))$ according to \eqref{eq:H_BC}
as well as
\begin{align}\notag%
	&\pmb{F}_H[G](\pxi,\peta) =\\ \notag
	&\frac{1}{4\ps} \int_{\peta_l(\pxi)}^{\peta} 
	\bigg(
		\big(
			\pa\big(z_{ij}(\pxi,\oeta),\zeta_{ij}(\pxi,\oeta)\big) + \mu_c
		\big) 
		\pG(\pxi,\oeta) \\\notag
		&+ \sum_{k=1}^{n} \Big(
			G_{ik}\big(
				\xi^{\xi\eta}_{ik}(\pxi,\oeta),\eta^{\xi\eta}_{ik}(\pxi,\oeta)
			\big) \,
			c^2_{kj}(\sigma_{ij}(\pxi,\oeta)) \\\notag
			&+ \int_{\zeta_{ij}(\pxi,\oeta)}^{z_{ij}(\pxi,\oeta)}
				c^3_{kj}\big(\sigma_{ij}(\pxi,\oeta),\ozeta\big) \\ \label{eq:int_operator_H} & \quad
				\cdot
				G_{ik}\big(
					 \xi_{ik}\big(z_{ij}(\pxi,\oeta),\ozeta\big),\eta_{ik}\big(z_{ij}(\pxi,\oeta),\ozeta\big)
				\big)
			\d\ozeta
		\Big)
	\bigg)
	\d\oeta.
\end{align}%

\subsection{Successive Approximation}\label{sec:sucappr}
To compute the solution of the kernel integral equations \eqref{eq:inteq1_norm} and \eqref{eq:inteq2_norm}, the method of successive approximations presented in \cite{Sm04} is extended to the considered multivariable case. However, there is not only a single integral equation, but $n^2$ integral equations for each $G$ and $H$, which are coupled. 

With the recursions
\begin{subequations}\label{eq:update_law_orig}
	\begin{align}
	\pG^{l+1} &= \widetilde{\pmb{F}}_G[G^l,\pH^l],  &&  l \in \mathbb{N}_0 \\
	\pH^{l+1} &= \widetilde{\pmb{F}}_H[G^l],  && G^0 = 0, \ H^0 = 0,
	\end{align}
\end{subequations}
implied by \eqref{eq:inteq1_norm} and \eqref{eq:inteq2_norm}, the corresponding fixpoints
\begin{align}
	\pG = \lim_{N\rightarrow \infty} \pG^N \quad \text{and}\quad \pH = \lim_{N\rightarrow \infty}\pH^N
\end{align}
are the solutions of the integral equations, if they exist. The latter amounts to proving the corresponding convergence, for which the equivalent representation
\begin{align} \label{eq:equiv_represenation}
\pG &= \sum\limits_{l=0}^\infty \Delta \pG^l \quad \text{and} \quad \pH = \sum\limits_{l=0}^\infty \Delta \pH^l 
\end{align}
with
\begin{subequations} \label{eq:update_law}
	\begin{align}\label{eq:update_law1}
	\Delta \pG^{l+1} &= \pmb{F}_G[\Delta G^l, \Delta \pH^l], \quad  &\Delta \pG^0 &= \pG_0(\peta) \\
	\Delta \pH^{l+1} &= \pmb{F}_{H}[\Delta G^l], &\Delta \pH^0 &= \pmb{H}_0(\pxi) \label{eq:update_law2}
	\end{align}
\end{subequations}
is considered.

As the integral operators \eqref{eq:int_operatorG} and \eqref{eq:int_operator_H} contain sums over the kernel elements that are defined in different coordinate systems, a growth assumption is needed, which is independent from the coordinate systems of the kernel elements. For this purpose, assume
\begin{subequations}\label{eq:grow}
	\begin{align}\label{eq:grow_G}
		|\Delta\pG^l(\xi_{ij}(z,\zeta),\eta_{ij}(z,\zeta))| \leq \frac{M^{l+1}}{l!}(z-\gamma \zeta)^l \\
		|\Delta\pH^l(\xi_{ij}(z,\zeta),\eta_{ij}(z,\zeta))| \leq \frac{M^{l+1}}{l!}(z-\gamma \zeta)^l  \label{eq:grow_H}
	\end{align}
\end{subequations}
with 
	\begin{align}\label{eq:gamma}
		\gamma \in 
		\left(\max\limits_{\lambda_i<\lambda_j}
		\sqrt{\frac{\lambda_i(\pmb{z}_{\underline{\Delta}})}{\lambda_j(\pmb{z}_{\underline{\Delta}})}},\ 1\right).
	\end{align}
Thereby,
\begin{equation}
	\pmb{z}_{\underline{\Delta}} = \operatorname{argmin}(|\lambda_i(z)-  \lambda_j(z)|)
\end{equation}
is the point of minimal difference between the diffusion coefficients $\lambda_i$ and $\lambda_j$. 

The \emph{growth assumptions} \eqref{eq:grow} can be inserted in \eqref{eq:equiv_represenation} and expressed in the canonical coordinates with \eqref{eq:inverse_one_step_trafo}.
This gives
\begin{subequations}\label{eq:boundedness}
\begin{align}
	|\pG(\pxi,\peta)|  &\leq \sum_{l=0}^{\infty} \frac{M^{l+1}}{l!}(z_{ij}(\pxi,\peta)-\gamma\zeta_{ij}(\pxi,\peta))^l  \\
	&=  M\e^{M(z_{ij}(\pxi,\peta)-\gamma\zeta_{ij}(\pxi,\peta))}
	= M\e^{M(z-\gamma\zeta)}
	 \notag \\
	|\pH(\pxi,\peta)|  &\leq \sum_{l=0}^{\infty} \frac{M^{l+1}}{l!}(z_{ij}(\pxi,\peta)-\gamma\zeta_{ij}(\pxi,\peta))^l  \\
	&=  M\e^{M(z_{ij}(\pxi,\peta)-\gamma\zeta_{ij}(\pxi,\peta))}= M\e^{M(z-\gamma\zeta)} \notag
\end{align}
\end{subequations}
showing that the series \eqref{eq:equiv_represenation}
converge absolutely and uniformly. 
Hence, the piecewise continuous solution of the kernel integral equations \eqref{eq:inteq1_norm} and \eqref{eq:inteq2_norm} can be computed via \eqref{eq:update_law_orig}, if
\eqref{eq:grow} is verified for all elements $\Delta \pG^l$ and $\Delta \pH^l$. 

Considering the initial values $\Delta\pG^0$ and $\Delta \pH^0$ given by \eqref{eq:G0} and \eqref{eq:H0}, it can easily be seen that \eqref{eq:grow} is valid for $l=0$, if the coefficient functions $\pmb{c}^{1}$, $\pmb{c}^{5}$, $\pmb{c}^{9}$ and $\pmb{c}^{10}$, as well as the degree of freedom $\pmb{g}_f$ are bounded. The definitions of the coefficient functions directly show that this is always true, if the system parameters $\Lambda$, $A$, $A_0$ and $F$  in \eqref{drs} are bounded.
What remains to be proven is that the integral operators \eqref{eq:int_operatorG} and \eqref{eq:int_operator_H} used for the update law \eqref{eq:update_law} preserve the growth assumption \eqref{eq:grow} for all elements $\Delta \pG^{l}$ and $\Delta \pH^{l}$, $l > 0$.

\subsubsection{Preliminaries}
To prove the validity of the growth assumptions \eqref{eq:grow} for $l >0$, it is convenient to take some preliminary algebraic considerations first. 


The inverse transformation \eqref{eq:inverse_one_step_trafo} can be inserted in $z-\gamma \zeta$ and differentiated \wrt the canonical coordinates $(\pxi, \peta)$. To shorten the presentation, the points $(\pz,\pzeta)$ will be written instead of the functions $z_{ij}(\pxi,\peta)$ and $\zeta_{ij}(\pxi,\peta)$ in the following. This yields
\begin{align}\label{eq:dxi_z_m_gamma_z}
	\partial_{\pxi} (\pz-\gamma \pzeta) &= 
	\partial_{\pxi} \big(\phi_i^{-1}(\tfrac{1}{2}(\ps\pxi+\peta)+\frac{1}{2}(1-\ps)\phi_i(1))\big) \\\notag
	&\quad -\gamma \partial_{\pxi} \big(\phi_j^{-1}(\tfrac{1}{2}(\ps\pxi-\peta)+\tfrac{1}{2}(1-\ps)\phi_j(1))\big).
\end{align}
Considering the well-known result for the derivative of the inverse function of \eqref{phidef}
\begin{align}
	\frac{\d}{\d y} \phi_i^{-1} (y)\big|_{y=\phi_i(z)} =  \frac{1}{\phi_i'(z)} = \sqrt{\lambda_i(z)},\quad  i=1,2,\ldots,n
\end{align}
and inserting \eqref{eq:inverse_one_step_trafo} leads to the result
\begin{align}\label{eq:dxzmgz>0}
	\partial_{\pxi} (\pz-\gamma \pzeta) = \tfrac{1}{2}\ps \Big(\sqrt{\lambda_{i\vphantom{j}}(\pz)} - \gamma \sqrt{\lambda_j(\pzeta)}\: \Big). 
\end{align}
With the definitions \eqref{shdef} and \eqref{eq:gamma} it is easy to see that \eqref{eq:dxzmgz>0} is strictly positive. Repeating this procedure for the derivative \wrt $\peta$ shows that in summary
\begin{subequations} \label{eq:dxi_deta>0}
	\begin{align}
		\partial_{\pxi} (\pz-\gamma \pzeta) > 0  \\
		\partial_{\peta} (\pz-\gamma \pzeta) > 0 
	\end{align}
\end{subequations}
hold. 
In view of the spatial domain $0\leq \zeta \leq z \leq 1$ in the original coordinates and \eqref{eq:gamma}, the inequality
\begin{align}\label{eq:ineq_gamma}
	0 \leq z-\gamma\zeta \leq 1
\end{align}
can be deduced. With this and \eqref{eq:dxi_deta>0}, the relations
\begin{subequations}\label{eq:relations}
	\begin{align}
		0 \leq \label{eq:relation1} \frac{z_{ij}(\xi^{\ast},\peta)-\gamma\zeta_{ij}(\xi^{\ast},\peta)}{z_{ij}(\pxi,\peta)-\gamma\zeta_{ij}(\pxi,\peta)} \leq 1\\ \label{eq:relation2}
		0 \leq \frac{z_{ij}(\pxi,\eta^{\ast})-\gamma\zeta_{ij}(\pxi,\eta^{\ast})}{z_{ij}(\pxi,\peta)-\gamma\zeta_{ij}(\pxi,\peta)} \leq 1
	\end{align}
\end{subequations}
are satisfied for all points $\xi^{\ast} \leq \pxi$ and $\eta^{\ast} \leq \peta$.

Moreover, calculating the derivative of $(\pz-\gamma\pzeta)^{l+1}$ \wrt $\pxi$ and using \eqref{eq:dxzmgz>0} results in
	\begin{align} \label{eq:eq_to_integrate}
		\partial_{\pxi}(\pz-\gamma \pzeta)^{l+1} 
		&= (l+1)(\pz-\gamma \pzeta)^{l} \partial_{\pxi}(\pz-\gamma \pzeta) \\ \notag
		&= (l+1)(\pz-\gamma \pzeta)^{l} \tfrac{\ps}{2} \Big(\sqrt{\lambda_{i\vphantom{j}}(\pz)} - \gamma \sqrt{\lambda_j(\pzeta)}\:\Big) \\ \notag
		&\hspace{-0.7cm} \geq (l+1)(\pz-\gamma \pzeta)^{l} \tfrac{1}{2} 
		\underbrace{
			\ps \Big(\sqrt{\lambda_i(\pz_{\underline{\Delta}})} - \gamma \sqrt{\lambda_j(\pz_{\underline{\Delta}})}\:\Big)
		}
		_{\big|\sqrt{\lambda_i(\pz_{\underline{\Delta}})} - \gamma \sqrt{\lambda_j(\pz_{\underline{\Delta}})}\big|}.
	\end{align}
At this point, the importance of the choice of $\gamma$ according to \eqref{eq:gamma} becomes visible. It ensures that the right sides of \eqref{eq:dxzmgz>0} and \eqref{eq:eq_to_integrate} are positive which is the basis for \eqref{eq:relations} and the following bound estimations.

Solving \eqref{eq:eq_to_integrate} for $(\pz - \gamma \pzeta) ^l$ and integrating both sides \wrt $\pxi$ leads to the upper bound
\begin{subequations}
	\begin{multline} \label{eq:int_xibar}
		\int_{\pxi_l}^{\pxi}\big(z_{ij}(\oxi,\peta) - \gamma \zeta_{ij}(\oxi,\peta)\big)^l \d\oxi \\
		\leq%
		\frac{2 \Big[%
			\big(z_{ij}(\oxi,\peta) - \gamma \zeta_{ij}(\oxi,\peta) \big) ^{l+1}%
		\Big]_{\pxi_l}^{\pxi}%
		}{
			(l+1)\left|\sqrt{\lambda_i(\pz_{\underline{\Delta}})} - \gamma%
			\sqrt{\lambda_j(\pz_{\underline{\Delta}})}\right|
		}.
	\end{multline}
With an identical reasoning, it is possible to prove the inequality
\begin{multline}\label{eq:int_etabar}
	\int_{\peta_l}^{\peta}\big(z_{ij}(\pxi,\oeta) - \gamma \zeta_{ij}(\pxi,\oeta)\big)^l \d\oeta \\
	\leq 
	\frac{2 \Big[ 
		\big(z_{ij}(\pxi,\oeta) - \gamma \zeta_{ij}(\pxi,\oeta) \big) ^{l+1}
	\Big]_{\peta_l}^{\peta}
	}{(l+1)\left(\sqrt{\lambda_i(\pz_{\underline{\Sigma}})} + \gamma \sqrt{\lambda_j(\pz_{\underline{\Sigma}})}\:\right)}
\end{multline}
and similarly for $\lambda_i \geq \lambda_j$
\begin{multline}\label{eq:int_etabarbar}
	\int_{0}^{\peta}\big(z_{ij}(\oeta,\oeta) - \gamma \zeta_{ij}(\oeta,\oeta)\big)^l \d\oeta \\
	\leq 
	\frac{\Big[ 
		\big(z_{ij}(\oeta,\oeta) - \gamma \zeta_{ij}(\oeta,\oeta) \big) ^{l+1}
	\Big]_{0}^{\peta}
	}{(l+1)\sqrt{\lambda_i(\pz_{\underline{\Sigma}})}},
\end{multline}
\end{subequations}
where 
\begin{align}
	\pz_{\underline{\Sigma}} = \operatorname{argmin}(\lambda_i(z)+\lambda_j(z))
\end{align}
is the point for which the sum of the diffusion coefficients $\lambda_i$ and $\lambda_j$ is minimal.

By making use of \eqref{eq:one_step_trafo} and \eqref{eq:inverse_one_step_trafo}, it can be verified that
\begin{subequations}\label{eq:coordsik}%
\begin{align}%
	& z_{ik}\big(\xi_{ik}(z_{ij}(\oeta,\oeta),\ozeta),\eta_{ik}(z_{ij}(\oeta,\oeta),\ozeta)\big) = z_{ij}(\oeta,\oeta) \\
	&\zeta_{ik}\big(\xi_{ik}(z_{ij}(\oeta,\oeta),\ozeta),\eta_{ik}(z_{ij}(\oeta,\oeta),\ozeta)\big) = \ozeta 
\end{align}%
\end{subequations}%
is satisfied, which will significantly simplify the evaluation of the sums in the following steps.
\subsubsection{Proof of the growth assumptions for $l>0$}
Now the growth assumptions \eqref{eq:grow} can be inserted into the fixpoint iteration \eqref{eq:update_law} with the operators \eqref{eq:int_operatorG} and \eqref{eq:int_operator_H}, to prove its validity via induction. 



With 
\eqref{eq:grow}, \eqref{eq:coordsik}, the 
triangle inequality and factoring $\frac{M^{l+1}}{l!}$ out, the recursion
\eqref{eq:update_law1} and \eqref{eq:int_operatorG} give
\begin{align} \notag 
	&|\Delta \pG^{l+1}(\pxi,\peta)| \leq
	\frac{M^{l+1}}{l!}
	\Big(
	\int_{\pxi_l}^{\pxi} \Big(z_{ij}(\oxi,\peta)-\gamma\zeta_{ij}(\oxi,\peta)\Big)^l\d\oxi \\\label{eq:dG1}
    &+ \Big[ 
	\int_0^{\peta}
	\Big(
		 (2+c^4_{j})
		 \big(z_{ij}(\oeta,\oeta)-\gamma\zeta_{ij}(\oeta,\oeta)\big)^l \\\notag
		 &+ \int_0^{z_{ij}(\oeta,\oeta)}
		 \sum_{k=1}^{n}
		 \big(
			 z_{ij}(\oeta,\oeta)-\gamma \ozeta
		 \big)^l
		 c^6_{kj}(\ozeta)  \d\ozeta
	\Big)
	\d\oeta
	\Big]\rightcond{\lambda_i&\geq\!\lambda_j\\ j&>\!m\\\peta &>\!0}
	\Big).
\end{align}

Due to the fact that the growth assumptions \eqref{eq:grow} were formulated in the original coordinates, which led to the simplification \eqref{eq:coordsik}, the term $(z_{ij}(\oeta,\oeta)-\gamma \ozeta)^l$
can now be factored out of the sum, \ie
\begin{align}\label{eq:sum_simp}
	\sum_{k=1}^{n} \notag
		( z_{ij}(\oeta,\oeta)-\gamma \ozeta
		)^l
		c^6_{kj}(\ozeta) &= ( z_{ij}(\oeta,\oeta)-\gamma \ozeta
		)^l 
	\sum_{k=1}^{n}	c^6_{kj}(\ozeta) \\  
	& \leq ( z_{ij}(\oeta,\oeta)-\gamma \ozeta
			)^l \, n\, \bar{c}
\end{align}
is obtained, in which
$c^6_{kj}(\ozeta) \leq \bar{c}$, $k,j = 1,2,\ldots,n$, $\ozeta \in [0,1]$, follows from the boundedness of the system parameters.


With \eqref{eq:sum_simp}, the remaining inner integral in \eqref{eq:dG1} can be evaluated as
\begin{align}\label{eq:inner_int}
	& \int_0^{z_{ij}(\oeta,\oeta)}
	\Big(z_{ij}(\oeta,\oeta)-\gamma\ozeta\Big)^l \d\ozeta  \\\notag&
	= -\frac{1}{\gamma (l+1)}\Big[\Big(z_{ij}(\oeta,\oeta)-\gamma\ozeta\Big)^{l+1}\Big]_0^{z_{ij}(\oeta,\oeta)} \\\notag&
	= \frac{1}{\gamma (l+1)}\Big(z_{ij}(\oeta,\oeta)^{l+1} - \big(z_{ij}(\oeta,\oeta) - \gamma z_{ij}(\oeta,\oeta)\big)^{l+1}\Big) \\\notag
	&= \frac{1}{\gamma(l+1)} \Big(z_{ij}(\oeta,\oeta)-\gamma\zeta_{ij}(\oeta,\oeta)\Big)^{l+1} \\\notag
	&
	\cdot\left(\!\!\!
	\left(\frac{z_{ij}(\oeta,\oeta)}{z_{ij}(\oeta,\oeta)-\gamma\zeta_{ij}(\oeta,\oeta)}\right)^{\:\: \mathclap{l+1}}
	\!\!\! - \!\left(
	\frac{z_{ij}(\oeta,\oeta)-\gamma z_{ij}(\oeta,\oeta)
	}{z_{ij}(\oeta,\oeta)-\gamma\zeta_{ij}(\oeta,\oeta)}
		\right)^{\!\! l+1}
	\right).
\end{align}
As the square bracket in \eqref{eq:dG1} only exists for $\lambda_i\geq \lambda_j$, the inverse transformation \eqref{eq:inverse_one_step_trafo} has to be considered for $\ps=1$ in \eqref{eq:inner_int} (see \eqref{shdef}). Consequently, inserting $\pxi=\peta=\oeta$ in \eqref{eq:zeta_inv} results in $\zeta_{ij}(\oeta,\oeta) = 0$. Thus, \eqref{eq:inner_int} can be simplified to
\begin{align} \notag
	& \int_0^{z_{ij}(\oeta,\oeta)}
	\Big(z_{ij}(\oeta,\oeta)-\gamma\ozeta\Big)^l \d\ozeta  \\\notag
	&\quad \leq \frac{1}{\gamma(l+1)} \big(z_{ij}(\oeta,\oeta)-\gamma\zeta_{ij}(\oeta,\oeta)\big)^{l+1}
	\big( 1-(1-\gamma)^{l+1} \big) \\\notag
	&\quad \leq \frac{1}{\gamma(l+1)} \big(z_{ij}(\oeta,\oeta)-\gamma\zeta_{ij}(\oeta,\oeta)\big)^{l+1} \\\label{eq:inner_int_simp}
	&\quad \leq \frac{1}{\gamma} \big(z_{ij}(\oeta,\oeta)-\gamma\zeta_{ij}(\oeta,\oeta)\big)^l,
\end{align} 
because of \eqref{eq:gamma}, \eqref{eq:ineq_gamma} and $l>0$.
The first appearance of $\zeta_{ij}(\oeta,\oeta)$ is kept deliberately, because it removes the necessity of case distinctions in the sequel. 
Inserting \eqref{eq:inner_int_simp} in \eqref{eq:dG1} yields
\begin{align} \notag 
	&|\Delta \pG^{l+1}(\pxi,\peta)| \leq
	\frac{M^{l+1}}{l!}
	\Big(
	\int_{\pxi_l}^{\pxi} \big(z_{ij}(\oxi,\peta)-\gamma\zeta_{ij}(\oxi,\peta)\big)^l\d\oxi \\
    &+ \Big[ 
    \underbrace{(2+c^4_{j} + \frac{n\:\bar{c}}{\gamma})}_{\bar{\bar{c}}}
	\int_0^{\peta}
	 \big(z_{ij}(\oeta,\oeta)-\gamma\zeta_{ij}(\oeta,\oeta)\big)^l
	 \d\oeta
	\Big]\rightcond{\lambda_i&\geq\!\lambda_j\\ j&>\!m\\\peta &>\!0} \Big).
	\label{eq:Glp1_mid}
\end{align}
To evaluate the integrals in \eqref{eq:Glp1_mid}, insert \eqref{eq:int_xibar} and \eqref{eq:int_etabarbar}. After factoring out $\frac{1}{l+1}(z_{ij}(\pxi,\peta)-\gamma\zeta_{ij}(\pxi,\peta))^{l+1}$, the relation
\begin{flalign}%
	\notag
	&|\Delta \pG^{l+1}(\pxi,\peta)| \leq \frac{M^{l+1}}{l!(l+1)}%
	\big(z_{ij}(\pxi,\peta)-\gamma\zeta_{ij}(\pxi,\peta)\big)^{l+1}\hspace{-1mm}\\\notag
	& \bigg(
		\frac{2}{|\sqrt{\lambda_i(\pz_{\underline{\Delta}})}\!-\!\gamma\sqrt{\lambda_j(\pz_{\underline{\Delta}})}|}
		\Big(\!
			1 \!-\! \Big(\frac{z_{ij}(\pxi_l,\peta)\!-\!\gamma\zeta_{ij}(\pxi_l,\peta)}{z_{ij}(\pxi,\peta)\!-\!\gamma\zeta_{ij}(\pxi,\peta)}\Big)^{l+1}
		\Big) \\\notag
		&
		+ \Big[ 
			\frac{\bar{\bar{c}}}{\sqrt{\lambda_i(\pz_{\underline{\Sigma}})}}
			\Big(
				\Big(
					\frac{z_{ij}(\peta,\peta) - \gamma \zeta_{ij}(\peta,\peta)}
					{z_{ij}(\pxi,\peta)-\gamma\zeta_{ij}(\pxi,\peta)}
				\Big)^{l+1} \\
		&\hspace{1.8cm}
				-\Big(
					\underbrace{\frac{z_{ij}(0,0) - \gamma \zeta_{ij}(0,0)}
					{z_{ij}(\pxi,\peta)-\gamma\zeta_{ij}(\pxi,\peta)}}
					_{=0}
				\Big)^{l+1}
			\Big)
		\Big]\rightcond{\lambda_i&\geq\!\lambda_j\\ j&>\!m\\\peta &>\!0}
	\bigg) \hspace{-1cm} \label{eq:dG1_2}
\end{flalign}
is obtained.
In order to prove convergence, it is necessary that the absolute values of all terms with the exponent $l+1$ are lower or equal to one. Observe that $\peta \leq \pxi$ in the domain $\mathcal{D}_{ij}$ (see Figure \ref{Fig1}) and obviously $\pxi_l \leq \pxi$ as it is the left boundary of the domain. Thus, \eqref{eq:relation1} holds for the fractions in \eqref{eq:dG1_2}, why it can be further simplified to
\begin{align} \notag
	|\Delta \pG^{l+1}(\pxi,\peta)| &\leq \frac{M^{l+1}}{(l+1)!}   \big(z_{ij}(\pxi,\peta)-\gamma\zeta_{ij}(\pxi,\peta)\big)^{l+1} \\ \notag
	&\mathpushright[0.72]{
	\cdot\underbrace{\bigg(
		\frac{2}{\sqrt{\lambda_i(\pz_{\underline{\Delta}})}-\gamma\sqrt{\lambda_j(\pz_{\underline{\Delta}})}} 
		+\Big[
			\frac{\bar{\bar{c}}}{\sqrt{\lambda_i(\pz_{\underline{\Sigma}})}}
		\Big]\rightcond{\lambda_i&\geq\!\lambda_j\\ j&>\!m\\\peta &>\!0}
	\bigg)}_{\leq \tilde{c}}
	} \\ 
	& \leq \frac{M^{l+1}\,\tilde{c}}{(l+1)!}   \big(z_{ij}(\pxi,\peta)-\gamma\zeta_{ij}(\pxi,\peta)\big)^{\mathrlap{l+1}},
\end{align}
where $\tilde{c}$ is a sufficiently large but finite number. Now an
\begin{align}\label{eq:M1}
M \geq \tilde{c}
\end{align}
can always be found such that
\begin{align}\label{eq:grow_prove_G}
	|\Delta \pG^{l+1}(\pxi,\peta)| \leq \frac{M^{l+2}}{(l+1)!}   \big(z_{ij}(\pxi,\peta)-\gamma\zeta_{ij}(\pxi,\peta)\big)^{l+1}
\end{align}
is satisfied.
Comparing \eqref{eq:grow_prove_G} with \eqref{eq:grow_G} shows the validity of the growth assumption \eqref{eq:grow_G} for all $l>0$ by induction.

Now \eqref{eq:grow_G} is inserted into the second update law \eqref{eq:update_law2} with the integral operator \eqref{eq:int_operator_H}. 
After applying the simplification \eqref{eq:coordsik}, factoring out $\frac{M^{l+1}}{l!}$ and introducing the upper bound $\breve{c}$ of all appearing coefficients, the result
\begin{align}\notag
	|\Delta \pH^{l+1}(\pxi,\peta)| & \leq \frac{M^{l+1}\breve{c}}{l!}
		\int_{\peta_l}^{\peta}
		\bigg(
			\big(z_{ij}(\pxi,\oeta)-\gamma\zeta_{ij}(\pxi,\oeta)\big)^{l} \\ \notag
			& +
			\sum_{k=1}^n \Big(
				\big(z_{ij}(\pxi,\oeta)-\gamma\zeta_{ij}(\pxi,\oeta)\big)^l \\  \label{eq:dH_mid}
				&+\int_{\zeta_{ij}(\pxi,\oeta)}^{z_{ij}(\pxi,\oeta)}
					\big(z_{ij}(\pxi,\oeta)-\gamma\ozeta\big)^l
				\d\ozeta
			\Big)
		\bigg)
		\d\oeta
\end{align}
follows. As all terms in the sum are independent of $k$, the sum can be replaced by a multiplication with $n$. Moreover, the inner integral can be evaluated in a similar way as in \eqref{eq:inner_int}, which shows that it can be replaced by $\frac{1}{\gamma}(z_{ij}(\pxi,\oeta)-\gamma\zeta_{ij}(\pxi,\oeta))^l$ in 
\eqref{eq:dH_mid}. In summary, the estimate
\begin{align}
	&|\Delta \pH^{l+1}(\pxi,\peta)| \\\notag
	&\mathpushright[0.9]{\leq{} \frac{M^{l+1}\breve{c}}{l!} 
	\int\limits_{\peta_l}^{\peta}
		(2\!+\!n\! +\!\tfrac{n}{\gamma})\big(z_{ij}(\pxi,\oeta)-\gamma\zeta_{ij}(\pxi,\oeta)\big)^{l}
		\d\oeta}
\end{align}
is obtained. For the evaluation of the integral, the relation \eqref{eq:int_etabar} is applied. Then $(z_{ij}(\pxi,\peta)-\gamma\zeta_{ij}(\pxi,\peta))^{l+1}$ is factored out after inserting the limits of the integral, which leads to
\begin{align}
	&|\Delta \pH^{l+1}(\pxi,\peta)| \leq \frac{M^{l+1}\breve{c}}{l!(l+1)} \label{eq:H_lp1}
	\big(z_{ij}(\pxi,\peta) \!- \! \gamma\zeta_{ij}(\pxi,\peta)\big)^{l+1} \\\notag
	& \frac{2(2+n+\tfrac{n}{\gamma})}{\sqrt{\lambda_i(\pz_{\underline{\Sigma}})}+\gamma\sqrt{\lambda_j(\pz_{\underline{\Sigma}})}} 
	\Big(
		1 \!-\! \Big(\frac{z_{ij}(\pxi,\peta_l)-\gamma\zeta_{ij}(\pxi,\peta_l)}{z_{ij}(\pxi,\peta)-\gamma\zeta_{ij}(\pxi,\peta)}\Big)^{l+1}	
	\Big).
\end{align}
It is clear that $\peta_l \leq \peta$, as it is the lower boundary of the domain $\mathcal{D}_{ij}$. Thus, the relation \eqref{eq:relation2} can be used for further simplifying \eqref{eq:H_lp1} to 
\begin{align}\label{eq:grow_prove_H}
	&|\Delta \pH^{l+1}(\pxi,\peta)| \leq \frac{M^{l+2}}{(l+1)!}
	\big(z_{ij}(\pxi,\peta) - \gamma\zeta_{ij}(\pxi,\peta)\big)^{l+1} 
\end{align}
for some 
\begin{align}\label{eq:M2}
	M \geq \frac{\breve{c}\,2\,(2+n+\tfrac{n}{\gamma})}{\sqrt{\lambda_i(\pz_{\underline{\Sigma}})}+\gamma\sqrt{\lambda_j(\pz_{\underline{\Sigma}})}}.
\end{align}
Observe that with the definition of $\gamma$ in \eqref{eq:gamma} and the requirements for $M$, given by \eqref{eq:M1} and \eqref{eq:M2}, an appropriate, bounded $M$ can always be found.
With \eqref{eq:grow_prove_H} the validity of the growth assumption \eqref{eq:grow_H} for all $l>0$ follows by induction.


As a result, the kernel integral equations \eqref{eq:inteq1_norm} and \eqref{eq:inteq2_norm} admit a piecewise continuous solution.
The next theorem states that the corresponding kernel $K(z,\zeta)$ is also a piecewise classical solution of \eqref{keq}.
\begin{thm}[Kernel equations]\label{thm:kequ}
	The kernel equations \eqref{keq} have a piecewise $C^2$-solution on the spatial domain $0 \leq \zeta \leq z \leq 1$.	
\end{thm}	
\begin{IEEEproof}
As a result of the previous paragraph, the kernel integral equations \eqref{eq:inteq1_norm} and \eqref{eq:inteq2_norm} admit a piecewise continuous solution. The regularity of the kernel obtained from the successive approximation remains to be verified. For this purpose, differentiate \eqref{eq:inteq1_norm} twice \wrt $\bm{\eta}$ as well as \eqref{eq:inteq2_norm} \wrt $\bm{\xi}$ and $\bm{\eta}$. By utilizing the piecewise continuity of $\bm{G}(\bm{\xi},\bm{\eta})$ and $\bm{H}(\bm{\xi},\bm{\eta})$ as well as the assumed regularity of $\Lambda$, $A$ and $F$, this yields piecewise continuous functions. By direct substitution of the kernel integral equations in the respective kernel equations it follows that \eqref{canonkeq} and consequently \eqref{keq} have a piecewise classical solution, which proves the theorem.
\end{IEEEproof}
\section{Closed-loop stability}\label{sec:clstab}
In order to prove the stability of the closed-loop system resulting from applying \eqref{sfeed} to \eqref{drs}, the bounded invertibility of the backstepping transformation \eqref{btrafo} has to be verified. For this, the \emph{inverse backstepping transformation}
\begin{equation}\label{btrafoinv}
x(z,t) = \tilde{x}(z,t) + \int_0^zL(z,\zeta)\tilde{x}(\zeta,t)\d \zeta
\end{equation}
with the integral kernel $L(z,\zeta) \in \mathbb{R}^{n \times n}$ is introduced. In the proof of the next theorem it is demonstrated that the transformation \eqref{btrafoinv} mapping the target system \eqref{drst} back into the original coordinates exists. Thereby, the integral kernel $L(z,\zeta)$ is a piecewise $C^2$-function implying the boundedness of \eqref{btrafoinv}. Hence, the exponential stability of \eqref{drst} implies the same in the original coordinates, which is the result of the following theorem.
\begin{thm}[Closed-loop stability]\label{thm:clstab}
Assume that $\mu_c > \mu_{\text{max}}$ (see Theorem \ref{thm:tstab}). Then, the closed-loop system \eqref{drs} and \eqref{sfeed} is exponentially stable in the $L_2$-norm $\|h\| = (\int_0^1\|h(z)\|^2_{\mathbb{C}^n}\d z)^{1/2}$, i.\:e.,
\begin{equation}
 \|x(t)\| \leq M\text{e}^{(\mu_{\text{max}}-\mu_c) t}\|x(0)\|, \quad t \geq 0
\end{equation}
for all $x(0) \in (H^2(0,1))^n$ compatible with the BCs of the closed-loop system \eqref{drs}, \eqref{sfeed} and an $M \geq 1$.
\end{thm}
\begin{IEEEproof}
The transformation \eqref{btrafoinv} maps \eqref{drst} into the original coordinates, if $L(z,\zeta)$ is the solution of the
\emph{inverse kernel equations}
\begin{subequations}\label{keqinv}
	\begin{alignat}{2}
	&&\Lambda(z) L_{zz}(z,\zeta) - (L(z,\zeta)\Lambda(\zeta))_{\zeta\zeta}  &= \mathcal{E}[L](z,\zeta)\label{kpdeinv}\\
	&& L(z,0)\Lambda(0)E_1 &= -\mathcal{F}[L](z)E_1 \label{kBC0inv}\\
	& \mathpushleft[0]{L_{\zeta}(z,0)\Lambda(0)E_2 + L(z,0)(\Lambda'(0)E_2 + \Lambda(0)E_2Q_0)}\label{kBC1inv}\\
	&& \quad  &= \big(A_0(z) + \mathcal{F}[L](z)\big)E_2\notag \\
	&\mathpushleft[0]{\Lambda(z)L'(z,z) \!+\! \Lambda(z)L_{z}(z,z) \!+\! L_{\zeta}(z,z)\Lambda(z) + L(z,z)\Lambda'(z)}\nonumber\\
	&&\quad  &= -(A(z) \!+ \!\mu_cI)\label{kBC2inv}\\
	&& L(z,z)\Lambda(z) - \Lambda(z) L(z,z) &= 0\label{kBC3inv}\\  
	&& L(0,0) &= 0\label{kICinv}
	\end{alignat}
\end{subequations}
with \eqref{kpdeinv} defined on $0 < \zeta < z < 1$ and
\begin{subequations}
\begin{align}
\mathcal{E}[L](z,\zeta) &= -(A(z) + \mu_cI)L(z,\zeta) - F(z,\zeta)\nonumber\\
& \quad \, - \int_{\zeta}^{z}F(z,\zeta')L(\zeta',\zeta)\d\zeta'\\
\mathcal{F}[L](z) &= \widetilde{A}_0(z) + \int_0^zL(z,\zeta)\widetilde{A}_0(\zeta)\d \zeta.
\end{align}	
\end{subequations}
They follow from the same reasoning as in Section \ref{sec:keq}.
Note, that therein $\widetilde{A}_0(z)$ depends on the kernel $K(z,\zeta)$ (see \eqref{eq:Aij_tilde}) and thus is assumed to be known. 
{Considering the component form of \eqref{keqinv},
the BCs originating from \eqref{kBC0inv} and \eqref{kBC1inv} for $\lambda_i < \lambda_j$ are automatically fulfilled by $\widetilde{A}_0(z)$, which can be shown using the reciprocity relation between $K(z,\zeta)$ and $L(z,\zeta)$ (cf. \cite[Ch. 4.5]{Kr08}). The remaining BVP}
has the same form as
\eqref{keq}, 
so that there exists a piecewise $C^2$-solution $L(z,\zeta)$ by Theorem \ref{thm:kequ}. Hence, \eqref{btrafoinv} is bounded so that the proof of the theorem follows from Theorem \ref{thm:tstab} by utilizing standard arguments.
\end{IEEEproof}
\begin{figure}%
	{\pgfplotsset{%
			dotted/.style={line cap=round, dash pattern=on 1pt off 6pt, line width=2.5pt},%
			every axis/.append style={%
				ytick distance = 1,
				axis equal,
				extra y ticks={0},
				extra y tick labels = {\phantom{$-0.555$}},
				extra x ticks={0},
				extra x tick labels={\phantom{$\phi_1(1)+\phi_2(1)$}},
				every axis y label/.append style={yshift=-7mm},
				every axis x label/.append style={yshift=2mm},
			},%
			every axis title/.append style={yshift=-2mm},
			every axis post/.append style={
				width=0.39\linewidth,height=0.39\linewidth,
				every x tick label/.style={below=8pt,anchor=base},
			},%
		}%
		\hspace{-4.5mm}%
		\input{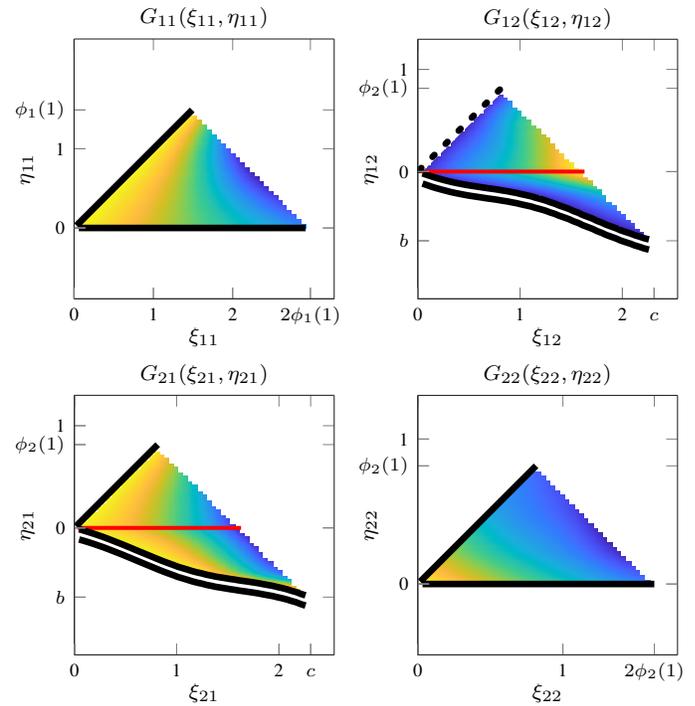}%
	}%
	\caption{Solution of the canonical kernel equations for $\mu_c=2$ with the color proportional to the surface hight ranging from blue to yellow. The solid lines represent the BCs, whereas the dashed line at $G_{12}(\xi_{12},\eta_{12})$ is the artificial BC with $g_f(\eta_{12})=0$ needed for well-posedness. The red lines in the spatial domains for the elements $G_{12}(\xi_{12},\eta_{12})$ and $G_{21}(\xi_{21},\eta_{21})$ are lines of separation. Furthermore, $b = \phi_2(1) - \phi_1(1)$ and $c=\phi_1(1)+\phi_2(1)$ specify the corresponding spatial domains.}
	\label{fig:G}
\end{figure}%

\section{Example}
Consider a system of two coupled PIDEs \eqref{drs} with the system parameters
\begin{subequations}
\begin{align}\notag
	\Lambda(z) &= 
	\begin{bmatrix}
		\tfrac{1}{2}-\tfrac{1}{4}\sin(2\pi z) & 0 \\
		0	& \tfrac{3}{2}+z^2 \cos(2\pi z) 
	\end{bmatrix}, \quad 	\Phi(z) = 0 \\
	A(z) &= \begin{bmatrix}
		1 & 1+z \\ \tfrac{1}{2}+z & 1
	\end{bmatrix}, \quad
	A_0(z) = \begin{bmatrix}
		z & 1-z \\ z & 1-z
	\end{bmatrix} \notag \\ \label{eq:ex}
	F(z,\zeta) &= \begin{bmatrix}
		\e^{z+\zeta} & \e^{z-\zeta} \\
		1-\e^{-(z-\zeta)} & \e^{-(z+\zeta)}
	\end{bmatrix}
\end{align}
and the differential operators
\begin{align}\label{eq:exBCleft}
	\theta_0[x(t)](0) &=  \begin{bmatrix}
	 0 & 0 \\0 & 1
	 \end{bmatrix} \dz x(0,t)  + \begin{bmatrix}
	  1 & 0 \\0 & -1
	  \end{bmatrix} x(0,t) \\ \label{eq:exBCright}
	 \theta_1[x(t)](1) &=
	  \begin{bmatrix}
	  1 & 0 \\ 0 & 0
	  \end{bmatrix} \dz x(1,t) + 
	  \begin{bmatrix}
	  0 & 1 \\
	  0 & 1
	  \end{bmatrix} x(1,t)
\end{align}
\end{subequations}
specifying the BCs, 
which is open-loop unstable. Note that $\lambda_2>\lambda_1$ is valid for \eqref{eq:ex}. 
The left BC \eqref{eq:exBCleft} specifies one Dirichlet and one Robin BC on the unactuated boundary. On the other side, \eqref{eq:exBCright} specifies Neumann actuation for the first state with a coupling of the second state which is actuated by a Dirichlet BC.
\begin{figure}
	{
		\pgfplotsset{
			dotted/.style={line cap=round, dash pattern=on 1pt off 6pt, line width=2.5pt},
			every axis/.append style={
				xtick distance=1,ytick distance=1,
				extra y ticks={0},
				extra y tick labels = {\phantom{$-0.555$}}, 
				every axis y label/.append style={yshift=-8mm},
				every axis x label/.append style={yshift=2mm}
			},
			every axis post/.append style={
				xmin=0, xmax = 1.23, ymin=-0.05,ymax=1.1,
				width=0.39\linewidth,height=0.39\linewidth,
			},
			every axis title/.append style={yshift=-2mm},
			every axis plot post/.append style={line cap=round}
		}
		\hspace{-5mm} 
		\input{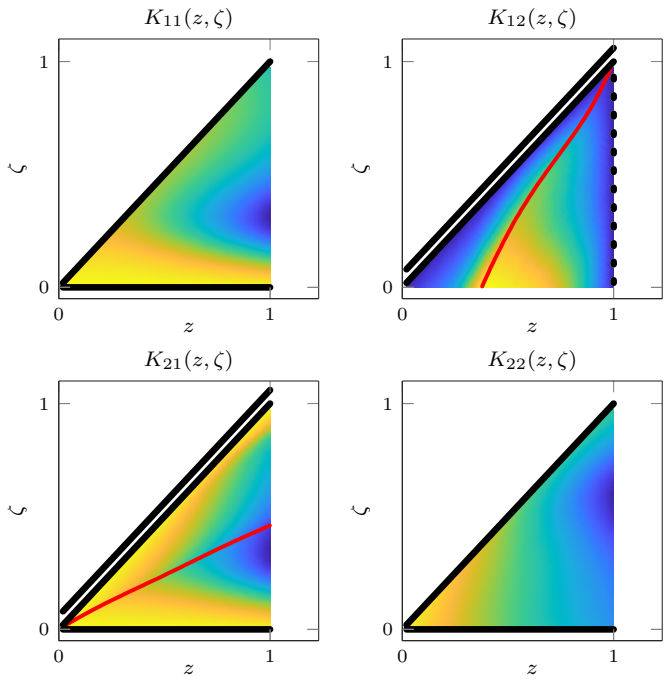}
	}
	\caption{Kernel elements in the original coordinates for $\mu_c=2$  with the color proportional to the surface hight ranging from blue to yellow. The solid lines represent the BCs, whereas the dashed line at $K_{12}(z,\zeta)$ is the artificial BC needed for well-posedness. The red lines in the spatial domains for the elements $K_{12}(z,\zeta)$ and $K_{21}(z,\zeta)$ are lines of separation.}
	\label{fig:K}
\end{figure}%
The controller is parametrized by $\mu_c = 2$ and the degree of freedom for specifying the artificial BC is chosen as $g_f(\eta_{12}) = 0$ (see \eqref{wpbc}). For the target system \eqref{tpdes} decoupled BCs at both sides are imposed. In particular, the same BC as in \eqref{eq:exBCleft} and
\begin{align}
	\tilde{\theta}_1[\tilde{x}(t)](1) = 
	\begin{bmatrix}
		1 & 0\\ 0 & 0
	\end{bmatrix} \dz \tilde{x}(1,t) +
	\begin{bmatrix}
		0 & 0\\0 & 1
	\end{bmatrix} \tilde{x}(1,t)
\end{align}
hold.

In order to solve the kernel equations, the successive approximation was numerically implemented in \textsc{Matlab}. Thereby, the spatial variables $z$ and $\zeta$ were discretised using 51 grid points each. 
The resulting grids for the different $(\xi,\eta)$-coordinate systems were resampled to ensure a minimum point distance for numerical performance. The successive approximation is stopped as soon as $\max_{\xi,\eta,i,j} \max(|\Delta G_{ij}^l|,|\Delta H_{ij}^l|)< 10^{-3}$, which occurs after 11 iteration steps.%

Figure \ref{fig:G} shows the solutions $G_{ij}(\xi_{ij},\eta_{ij})$, $i,j=1,2$, of the canonical kernel equations \eqref{canonkeq} in the respective spatial domains. Those are triangular domains for the diagonal elements. For the off-diagonal elements, the lower boundary evolves in the fourth quadrant because of the mutually different diffusion coefficients and is a strictly monotonically decreasing non-straight line as the diffusion coefficients are spatially-varying. 
The BCs required for the mapping into the target system are marked by the solid lines, whereas an artificial BC has to be assigned for the element $G_{12}$ (dotted line). Since the spatial domains of the elements $G_{12}$ and $G_{21}$ also cover the fourth quadrant, these kernel elements are only defined piecewise. This is indicated by the line of separation (red line), which is given by $\eta_{12}=0$ and $\eta_{21}=0$. On these lines, the obtained kernel elements are still continuous, but not differentiable.

In Figure \ref{fig:K}, the solution of the kernel equations $K_{ij}(z,\zeta)$, $i,j = 1,2$, in the original coordinates is depicted.
Here, the line of separation (red line) is no longer a straight line, but a strictly monotonically increasing curve due to the nonlinear change of coordinates.

\begin{figure}
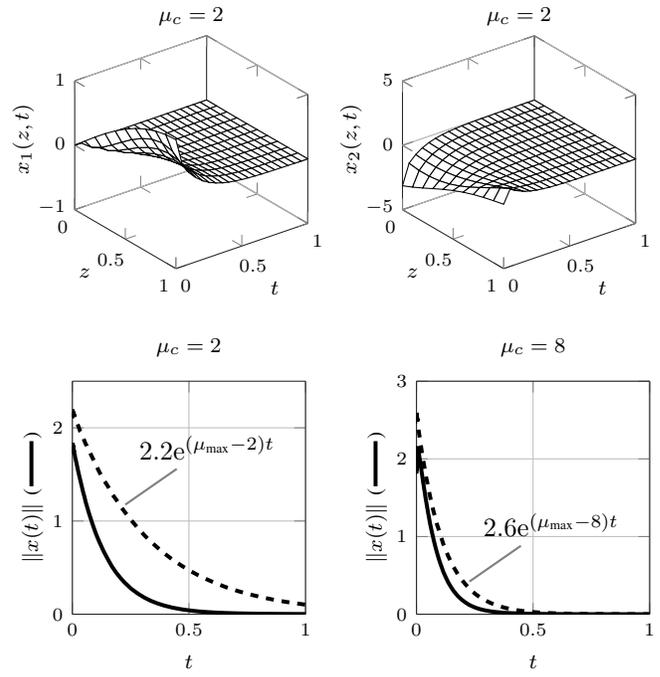

{\pgfplotsset{
	surf/.append style={shader = faceted,faceted color=black, fill=white},
	every axis plot post/.append style={shader=faceted, z buffer=auto},
	every axis/.append style={
		width=0.4\linewidth,
		height=0.4\linewidth,
		every axis y label/.append style={yshift=3mm, xshift=0mm},
		every axis x label/.append style={xshift=0mm,yshift=2mm},
		every axis z label/.append style={yshift=-2mm,xshift=0mm},
		extra z ticks={0},
		extra z tick labels = {\phantom{$-5$}},
		clip=false,
		execute at end axis={
			\node at (rel axis cs: -0.5,-0.5,0){{}};
		}
	},
	every axis title/.append style={yshift=-2mm},
	every axis post/.append style={
		xmajorgrids=false,
		ymajorgrids=false,
		zmajorgrids=true,
	}
}
\input{1_Verlauf_geregeltes_System.tex} 
}
\\
{
\newcommand{\dist}{0.5cm}
\pgfplotsset{
	every axis y label/.append style={yshift=-7mm},
	every axis x label/.append style={yshift=1mm},
	every axis/.append style={
		xtick distance = 0.5,
		execute at end axis={\hspace{\dist}}
	},
	every axis plot post/.append style={color=black},
}
\hspace*{-\dist}
\input{6_Norm_geregeltes_System}
}
\caption[Solution of the closed-loop system for $x_0(z)$]{Solution of the closed-loop system for $x_0(z) = \begin{bmatrix}
	\sin(\tfrac{\pi}{2}z) & \sin(\pi z) - \pi \cos(\tfrac{\pi}{2}z) 
\end{bmatrix}^{\top}$.
The first row depicts the closed-loop state profile for $\mu_c = 2$, whereas the plots in the second row contain the $L_2$-norms of the solution for $\mu_c=2$ and $\mu_c = 8$ with $\mu_\text{max} = -1.36$.}
\label{fig:x}
\end{figure}
For the simulation, the plant was discretised using a finite-element method with 102 grid points for each state.
To visualize the influence of the controller parameter $\mu_c$, Figure \ref{fig:x} shows the closed-loop state profiles for $\mu_c=2$ and the related $L_2$-norms for $\mu_c=2$ and $\mu_c = 8$. In both cases, the initial value is $x_0(z) = \begin{bmatrix}
	\sin(\tfrac{\pi}{2}z) & \sin(\pi z) - \pi \cos(\tfrac{\pi}{2}z) 
\end{bmatrix}^{\top}$, which is compatible to the BCs.
The result verifies that exponential stability is achieved and that a desired closed-loop stability margin can be assigned by $\mu_c$.

\section{Concluding remarks}
An interesting research topic is to take a convective coupling in the PIDEs and a  coupling at the unactuated boundary into account. As far as the solution of the kernel equations is concerned such an extension seems to be directly possible. However, since the aforementioned couplings have to be included in the target system, this yields, in general, a bidirectionally coupled system of parabolic PDEs. Hence, a different proof of well-posedness and stability is needed. For this, the corresponding results in \cite{Vaz16a} are of interest. In order to determine output feedback controllers the design of backstepping observers for the considered system class is currently under investigation. 

%
\IEEEpeerreviewmaketitle

\appendix\label{app:proofth1}
\emph{Proof of Theorem \ref{thm:tstab}}. 
Introduce the operators $\mathcal{A}_ih = \lambda_i\d^2_zh - \mu_ch$, $\mu_c > \mu_{\text{max}}$, $i = 1,2,\ldots,n$, with $D(\mathcal{A}_i) = \{h \in L_2(0,1) \;|\; h, \d_zh \; \text{abs. cont.}, \d_z^2h \in L_2(0,1), e_i^{\top}\theta{_0}[P^{\top}e_ih](0) = e_i^{\top}\tilde{\theta}_{1}[P^{\top}e_ih](1) = 0\}$ in the Hilbert space $H = L_2(0,1)$ with the usual inner product and $e_i \in \mathbb{R}^n$ denoting the $i$-th unit vector. 
Therein the largest eigenvalue $\mu_{\text{max}} \in \mathbb{R}$ of $\mathcal{A}_i + \mu_cI$, $i = 1,2,\ldots,n$, exists because $-(\mathcal{A}_i+ \mu_cI)$ are \emph{Sturm-Liouville operators} (see \cite{Del03}). 
With this, the operator $\mathcal{A}h = [\mathcal{A}_1h_1 \;\; \ldots \;\;\mathcal{A}_nh_n]^{\top}$, $D(\mathcal{A}) = D(\mathcal{A}_1) \oplus \ldots \oplus D(\mathcal{A}_n)$ can be defined. 
Furthermore, let $\Delta h = -\widetilde{A}^*_0(R_1\d_zh(0) + R_2h(0))$ with $D(\Delta) = D(\mathcal{A})$. Based on these preparations the target system \eqref{drst2} can be formulated in form of the abstract initial value problem (IVP) $\dot{\tilde{x}}(t) = (\mathcal{A} + \Delta)\,\tilde{x}(t)$, $t > 0$, $\tilde{x}(0) = \tilde{x}_0 \in D(\mathcal{A})$, on the state space $X = (L_2(0,1))^n$ with the $L_2$-inner product $\langle\cdot,\cdot\rangle$. 
Therein, $\mathcal{A}$ is the infinitesimal generator of an \emph{analytic $C_0$-semigroup}, because the operators $\mathcal{A}_i$ are the generators of analytic $C_0$-semigroups. The latter property follows from the fact that $\mathcal{A}_i$ are \emph{Riesz-spectral operators} satisfying the sector condition of \cite[Ex. 2.18]{Cu95}. 
Moreover, the operator $\Delta$ is \emph{$\mathcal{A}$-bounded} (see \cite[Ch. IV, \S1, Sec. 2]{Kat95}) and thus also \emph{$\mathcal{A}$-compact} due to its finite-dimensional range  (see \cite[Rem. IV/1.13]{Kat95}). This implies an $\mathcal{A}$-bound equal to zero (see \cite[Lem. III/2.16]{Eng00}). 
Thus, the perturbed operator $\mathcal{A} + \Delta$ is the generator of an analytic $C_0$-semigroup by \cite[Th. III/2.10]{Eng00}. 
Hence, the IVP \eqref{drst2} is well-posed. 
The operators $\mathcal{A}_i$ have a compact resolvent, since they have no eigenvalue at 0 (see \cite[Th. 7.5.4]{Na82}). 
Thus, also $\mathcal{A}$ has a compact resolvent, so that the $\mathcal{A}$-boundedness of $\Delta$ implies that $\mathcal{A} + \Delta$ has a compact resolvent (see \cite[Th. IV/3.17]{Kat95}). 
Therefore, the spectrum of the operator $\mathcal{A} + \Delta$ is discrete (see \cite[Th. III/6.29]{Kat95}). 
The \emph{spectrum determined growth assumption (SDGA)} holds for $\mathcal{A} + \Delta$ because it is a generator of an analytic semigroup (see \cite{Tr75}). 
This and the triangular structure of $\mathcal{A} + \Delta$ yield the decay rate $\mu_c-\mu_{\text{max}}$ for the semigroup related to $\mathcal{A} + \Delta$ and thus the stability result of the theorem. \hfill $\blacksquare$


\bibliographystyle{IEEEtranS}
\bibliography{IEEEabrv,mybib}

\begin{IEEEbiography}[{\includegraphics[width=1in,height=1.25in,clip,keepaspectratio]{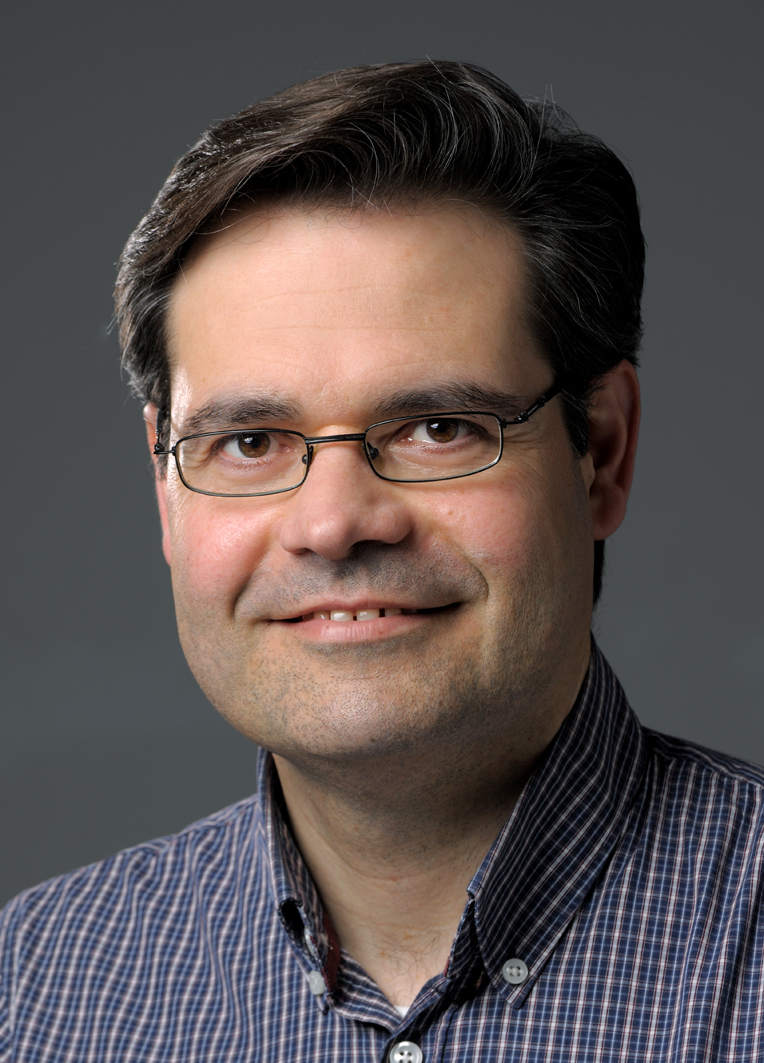}}]{Joachim Deutscher}
received the Dipl.-Ing. (FH) degree in electrical engineering from the University of Applied Sciences W\"urzburg-Schweinfurt-Aschaffenburg, Germany, in 1996, the Dipl.-Ing. Univ. degree in electrical engineering, the Dr.-Ing. and the Dr.-Ing. habil. degrees both in automatic control from the Friedrich-Alexander University Erlangen-Nuremberg (FAU), Germany, in 1999, 2003 and 2010, respectively.

From 2003--2010 he was a Senior Researcher at the Chair of Automatic Control (FAU), in 2011 he was appointed Associate Professor and since 2017 he is a Professor at the same university. Currently, he is Head of the Infinite-Dimensional Systems Group at the Chair of Automatic Control (FAU). 

His research interests include control of distributed-parameter systems and control theory for nonlinear lumped-parameter systems with applications in mechatronic systems. 

Dr. Deutscher has co-authored a book on state feedback control for linear lumped-parameter systems: Design of Observer-Based Compensators (Springer, 2009) and is author of the book: State Feedback Control of Distributed-Parameter Systems (in German) (Springer, 2012). 
\end{IEEEbiography}

\begin{IEEEbiography}[{\includegraphics[width=1in,height=1.25in,clip,keepaspectratio]{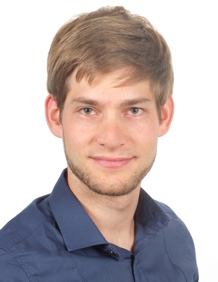}}]{Simon Kerschbaum} received the B.Sc. and M.Sc. degrees in electrical engineering from the Friedrich-Alexander University Erlangen-Nuremberg (FAU), Germany, in 2012 and 2014. Since then, he is a Ph.D. student at the Chair of Automatic Control (FAU) in the research group of Dr. Deutscher.
	
His research interests include backstepping methods for coupled parabolic systems and output regulation for distributed-parameter systems. 	
\end{IEEEbiography}
\vfill	
\end{document}